\documentclass{amsart}
\usepackage{amsmath,amssymb,amscd,amsthm}

\newtheorem{thm}{Theorem}[section]
\newtheorem{pro}[thm]{Proposition}
\newtheorem{lem}[thm]{Lemma}
\newtheorem{cor}[thm]{Corollary}
\newtheorem{dfn}[thm]{Definition}
\newtheorem{ex}[thm]{Example}
\newtheorem{rem}[thm]{Remark}

\begin{document}

\title{Weakly reflective submanifolds and austere submanifolds}

\author{Osamu Ikawa}
\address{Department of General Education,
Fukushima National College of Technology,
Iwaki, Fukushima, 970-8034 Japan}
\email{ikawa@fukusima-nct.ac.jp}

\author{Takashi Sakai}
\address{Graduate School of Science,
Osaka City University,
3-3-138 Sugimoto, Sumiyoshi-ku, Osaka-shi, Osaka, 558-8585 Japan}
\email{tsakai@sci.osaka-cu.ac.jp}

\author{Hiroyuki Tasaki}
\address{Graduate School of Pure and Applied Science,
University of Tsukuba,
Tsukuba, Ibaraki, 305-8571 Japan}
\email{tasaki@math.tsukuba.ac.jp}




\begin{abstract}
We introduce the notion of a weakly reflective submanifold,
which is an austere submanifold with a certain global condition,
and study its fundamental properties.
Using these, we determine weakly reflective orbits and austere orbits
of $s$-representations.
\end{abstract}

\maketitle

\section{Introduction}

Orbits of an $s$-representation,
that is a linear isotropy representation of a Riemannian symmetric pair,
are important examples of homogeneous submanifolds in the hypersphere
of a Euclidean space.
For example, a homogeneous isoparametric hypersurface in the hypersphere,
which many mathemtaticians have investigated,
can be obtained as a principal orbit of an $s$-representation
of a Riemannian symmetric pair of rank two.
The family of isoparametric hypersurfaces has a unique
minimal isoparametric hypersurface.
Furtheremore, typical examples of minimal submanifolds
in the hypersphere are given as orbits of $s$-representations.
Hirohashi-Song-Takagi-Tasaki \cite{HSTT00} showed that
there exists a unique minimal orbit in each strata of the stratification
of orbit types.
However, in general we can not explicitly point out which orbit
among each strata is a minimal submanifold.

Harvey-Lawson \cite{HL82} introduced the notion of an austere
submanifold, which is a minimal submanifold whose second fundamental form
has a certain symmetry.
They showed that one can construct a special Lagrangian cone,
therefore absolutely area-minimizing, in a complex Euclidean space
as the twisted normal bundle of an austere submanifold in a sphere
(see \cite{HL82}, \cite{BG2004}).
As we mentioned above, the complete list of
minimal orbits of $s$-representations in the hypersphere is unknown at the moment.
Therefore we first attempt to determine all austere orbits.
We give a necessary and sufficient condition for an orbit to be an austere
submanifold in the hypersphere in terms of the restricted root system of a Riemannian
symmetric pair.
By this criterion, we can determine all orbits which are austere submanifolds
in the hypersphere.
Since the definition is focused on a symmetry of its second fundamental form,
the notion of an austere submanifold is an infinitesimal property
of a submanifold.
However, we observe that some of austere orbits, which we classified,
have a certain global symmetry.
This symmetry is a globalization of the notion of an austere submanifold and
a weakened condition of a reflective submanifold.
Therefore we shall call them weakly reflective submanifolds,
and study some fundamental properties of them.
Finally we determine all weakly reflective orbits
of $s$-representations.

The organization of this paper is as follows.
In Section $2$,
we will give the definition of a weakly reflective submanifold
(Definition~\ref{dfn:weakly reflective submanifold}),
and recall some related notions.
We study their relationship and fundamental properties.
In Section $3$,
we summarize geometry of orbits
of $s$-representations of Riemannian symmetric pairs.
This will be a preliminary for the sections below.
In Section $4$,
we shall give the list of orbits of $s$-representations
which are weakly reflective submanifolds in the hypersphere
(Theorem~\ref{thm:weakly reflective orbit}).
We show that these orbits are weakly reflective submanifold
in the hypersphere there, however,
we will show that the list gives all weakly reflective orbits later.
In Section $5$,
we will give a criterion of austere orbits (Lemma~\ref{lem:A}),
and determine all orbits which are austere submanifolds
in the hypersphere (Theorem~\ref{thm:austere}).
Furthermore we show that austere orbits which are not enumerated
in the list of weakly reflective orbits are not weakly reflective submanifolds.
Then we will complete the proof of the list of weakly reflective orbits.
In Section $6$,
we will study relationships between weakly reflective submanifolds in a sphere
and those in Euclidean spaces or complex projective spaces.

The authors are profoundly grateful to Makoto Kimura and Osami Yasukura
for their helpful suggestion on Proposition \ref{pro:orbit of root}.
Before we wrote this paper, Kimura, Yasukura and the third named author
showed a previous version of Proposition \ref{pro:orbit of root}
which is unpublished, that is,
the orbit of the highest root of a compact Lie group
under the adjoint action is an austere submanifold in the hypersphere.
The authors would also like to thank Reiko Miyaoka for her valuable comments.
In fact, Proposition \ref{pro:principal orbit of cohomogeneity one actions}
was essentially suggested by her.

\section{Definitions and fundamental results}

We begin with recalling the definition of a reflective submanifold
given by Leung \cite{Leung73}.
Let $\tilde M$ be a complete Riemannian manifold.
A connected component of the fixed point set of an involutive isometry
of $\tilde M$ is called a {\it reflective submanifold}.
A reflective submanifold is a complete totally geodesic submanifold.
The involutive isometry which defines a reflective submanifold $M$
can be determined uniquely.
We call it the {\it reflection} of $M$ and denote by $\sigma_M$.
If $M$ is a reflective submanifold in $\tilde M$
and $\sigma_M$ is its reflection,
then for any normal vector $\xi \in T_x^\perp M$
$$
\sigma_M(x) = x, \qquad
(d\sigma_M)_x \xi = - \xi, \qquad
\sigma_M(M) = M
$$
hold.
Taking notice of these properties,
we define a weakly reflective submanifold
as follows.

\begin{dfn} \rm \label{dfn:weakly reflective submanifold}
Let $M$ be a submanifold of a Riemannian manifold $\tilde M$.
For each normal vector $\xi \in T_x^\perp M$ at each point $x \in M$,
if there exists an isometry $\sigma_\xi$ of $\tilde M$ which satisfies
$$
\sigma_\xi(x) = x, \qquad
(d\sigma_\xi)_x \xi = - \xi, \qquad
\sigma_\xi(M) = M,
$$
then we call $M$ a {\it weakly reflective submanifold}
and $\sigma_\xi$ a {\it reflection} of $M$ with respect to $\xi$.
\end{dfn}
In the case where $M$ is a hypersurface,
$\sigma_\xi$ is independent of the choice of $\xi$ at each point $x$.
In this paper mainly we deal with orbits of some isometric actions of
compact Lie groups.
We note that if $M$ is an extrinsic homogeneous submanifold in $\tilde M$,
that is an orbit of an isometric action of a Lie group on $\tilde M$,
then it suffices to show the condition to be a weakly reflective
submanifold at one point of $M$.

\begin{rem} \rm
For a reflective submanifold,
there exists a reflection which is independent of the choice of
a normal vector.
So it is clear that the definition of a weakly reflective submanifold
is a weakened condition of a reflective submanifold.
\end{rem}

\begin{ex} \label{ex:1}
$$
S^{n-1}(1) \times S^{n-1}(1) = \left\{ (x, y) \mid x, y \in S^{n-1}(1) \right\}
$$
is a weakly reflective submanifold
in $(2n-1)$-dimensional sphere $S^{2n-1}(\sqrt2)$ of radius $\sqrt 2$.
\end{ex}

\begin{proof}
Since $S^{n-1}(1) \times S^{n-1}(1)$ is a homogeneous submanifold
of $S^{2n-1}(\sqrt2)$,
it suffices to show the condition to be weakly reflective at
one point of $S^{n-1}(1) \times S^{n-1}(1)$.
The tangent space of $S^{n-1}(1) \times S^{n-1}(1)$ at
$$
x = (1, 0, \dots, 0, \stackrel{n+1\atop\smile}{1}, 0, \dots, 0)
\in S^{n-1}(1) \times S^{n-1}(1)
$$
is given by
$$
T_x(S^{n-1}(1) \times S^{n-1}(1))
= \{(0, x_2, \dots, x_n, 0, y_2, \dots, y_n)
\mid x_i, y_j \in \mathbf R\},
$$
and the normal space in $S^{2n-1}(\sqrt2)$ is
$$
T_x^\perp(S^{n-1}(1) \times S^{n-1}(1))
= \mathbf R(1, 0, \dots, 0, \stackrel{n+1\atop\smile}{-1}, 0, \dots, 0).
$$
Now we define an isometry $\sigma$ of $S^{2n-1}(\sqrt2)$ by
$$
\sigma(x_1, \dots, x_n, y_1, \dots, y_n)
= (y_1, \dots, y_n, x_1, \dots, x_n)
$$
for $(x_1, \dots, x_n, y_1, \dots, y_n) \in S^{2n-1}(\sqrt2)$.
Then
$$
\sigma(x) = x, \qquad
\sigma(S^{n-1}(1) \times S^{n-1}(1)) = S^{n-1}(1) \times S^{n-1}(1)
$$
and $d\sigma_x$ acts on $T_x^\perp(S^{n-1}(1) \times S^{n-1}(1))$ as $-\mbox{id}$.
Thus $S^{n-1}(1) \times S^{n-1}(1)$ is a weakly reflective submanifold
in $S^{2n-1}(\sqrt2)$.
\end{proof}

\begin{dfn} \rm
Let $M$ be a submanifold of a Riemannian manifold $\tilde M$.
We denote the shape operator of $M$ by $A$.
$M$ is called an {\it austere submanifold}
if for each normal vector $\xi \in T^\perp_x M$,
the set of eigenvalues of $A_\xi$ is invariant
(concerning multiplicities) under multiplication by $-1$.
It is obvious that an austere submanifold is a minimal submanifold.
\end{dfn}

The notion of an austere submanifold was first
given by Harvey-Lawson \cite{HL82}.

\begin{pro}\label{pro:weakly reflective is austere}
A weakly reflective submanifold is an austere submanifold.
\end{pro}

\begin{proof}
Let $M$ be a weakly reflective submanifold
in a Riemannian manifold $\tilde M$.
Then for each normal vector $\xi \in T_x^\perp M$,
there exists an isometry $\sigma_\xi$ of $\tilde M$ which satisfies
$$
\sigma_\xi(x) = x, \qquad
(d\sigma_\xi)_x \xi = - \xi, \qquad
\sigma_\xi(M) = M.
$$
For a normal vector $\xi \in T_x^\perp M$,
we denote by $A_\xi$ the shape operator of $M$ with respect to $\xi$
and by $h$ the second fundamental form of $M$.
For $X, Y \in T_xM$, we take vector fields $\tilde X$ and $\tilde Y$
defined on a neighborhood of $x$ in $\tilde M$ which are tangent to $M$
and $\tilde X_x = X$ and $\tilde Y_x = Y$.
Since $\sigma_\xi$ satisfies $\sigma_\xi(M) = M$,
vector fields $d\sigma_\xi\tilde X$ and $d\sigma_\xi\tilde Y$ are tangent to $M$.
Let $\bar\nabla$ denote the covariant derivative of $\tilde M$.
Then we have
\begin{eqnarray*}
h((d\sigma_\xi)_xX, (d\sigma_\xi)_xY)
&=&
(\bar\nabla_{d\sigma_\xi\tilde X}d\sigma_\xi\tilde Y)_x^\perp
= ((d\sigma_\xi)_x\bar\nabla_{\tilde X} \tilde Y)^\perp \\
&=&
(d\sigma_\xi)_x(\bar\nabla_{\tilde X} \tilde Y)^\perp
= (d\sigma_\xi)_x h(X, Y).
\end{eqnarray*}
From
\begin{eqnarray*}
\langle A_\xi(d\sigma_\xi)_xX, (d\sigma_\xi)_xY\rangle
&=& \langle h((d\sigma_\xi)_xX, (d\sigma_\xi)_xY), \xi\rangle \\
&=& \langle (d\sigma_\xi)_xh(X, Y), \xi\rangle
= \langle h(X, Y), (d\sigma_\xi)_x^{-1}\xi\rangle \\
&=& \langle h(X, Y), -\xi\rangle
= -\langle A_\xi X, Y\rangle,
\end{eqnarray*}
we have
$(d\sigma_\xi)_x^{-1}A_\xi(d\sigma_\xi)_x = -A_\xi$.
This implies that $(d\sigma_\xi)_x$ provides an isomorphism between
eigenspaces of $A_\xi$ for eigenvalues $\lambda$ and $-\lambda$.
Thus $M$ is an austere submanifold.
\end{proof}

In the rest of this section,
we shall study weakly reflective orbits of isometric actions of Lie groups
on Riemannian manifolds.
First we shall provide some preliminaries.
Let $G$ be a Lie group acting isometrically
on a Riemannian manifold $\tilde M$
and $G_x$ be the isotropy subgroup at $x$,
that is, $G_x = \{ g \in G \mid gx = x \}$.
Then the orbit $G(x)$ is diffeomorphic to the coset manifold $G/G_x$.
An orbit $G(x)$ is a {\it principal orbit} if, for any $y \in \tilde M$,
there exists $g \in G$ such that $G_x \subset gG_yg^{-1}$.
It is known that there exists a principal orbit.
The codimension of a principal orbit is called the {\it cohomogeneity}
of the action of $G$ on $\tilde M$.
An orbit which is not principal is called a {\it singular orbit}.
The differential of the action of $G_x$ defines a linear representation
of $G_x$ on $T_x \tilde M$ called the {\it linear isotropy representation}.
The tangent space $T_x(G(x))$ and the normal space $T_x^\perp(G(x))$
of $G(x)$ at $x$ are invariant subspaces of the linear isotropy representation.
The restriction of the linear isotropy representation to $T_x^\perp(G(x))$
is called the {\it slice representation} at $x$.

\begin{thm}{\rm (Slice representation theorem
\cite[Theorem 1.1]{Kollross2002},
\cite[Theorem 4.6]{PT87},
\cite[Proposition 5.4.7]{PT88}).}
The cohomogeneity of a slice representation equals the cohomogeneity of
the action of $G$ on $\tilde M$.
Moreover, $G(x)$ is a principal orbit
if and only if the slice representation at $x$ is trivial.
\end{thm}

\begin{pro}\label{pro:cohomogeneity one action}
Any singular orbit of a cohomogeneity one action
on a Riemannian manifold is a weakly reflective submanifold.
\end{pro}

\begin{proof}
Suppose that the isometric action of a Lie group $G$
on a Riemannian manifold $\tilde M$ is cohomogeneity one.
Let $G(x)$ be a singular orbit.

First we consider the case
where the codimension of $G(x)$ is equal or greater than $2$.
From the slice representation theorem,
the isotropy subgroup $G_x$ acts transitively
on the hypersphere in $T_x^\perp(G(x))$.
In particular, for any $\xi \in T_x^\perp(G(x))$ there exists $g \in G_x$
such that $dg_x(\xi) = -\xi$.
Therefore $g$ becomes a reflection of $G(x)$ at $x$ with respect to $\xi$.
Since $G(x)$ is a homogeneous submanifold,
$G(x)$ has a reflection with respect to any normal vector at any point.
Thus $G(x)$ is a weakly reflective submanifold in $\tilde M$.

When the codimension of $G(x)$ is $1$,
the slice representation at $x$ is not trivial.
Therefore for any $\xi \in T_x^\perp(G(x))$ there exists $g \in G_x$
such that $dg_x(\xi) = -\xi$.
Thus,
by the same discussion with above,
$G(x)$ is a weakly reflective submanifold in $\tilde M$.
\end{proof}

\begin{rem} \rm
Podest\'a \cite{Podesta97} proved that any singular orbit
of a cohomogeneity one action is an austere submanifold.
However, essentially he showed Proposition \ref{pro:cohomogeneity one action}.
\end{rem}

\begin{pro}\label{pro:principal orbit of cohomogeneity one actions}
Let $G$ be a connected Lie group acting isometrically on a complete,
connected Riemannian manifold $\tilde M$.
Suppose that the action of $G$ on $\tilde M$ is cohomogeneity one
with two singular orbits.
If there exists a principal orbit
which is a weakly reflective submanifold in $\tilde M$,
then it has a same distance from two singular orbits
and two singular orbits are isometric.
\end{pro}

\begin{proof}
Since there exist two singular orbits,
the orbit space $\tilde M/G$ is homeomorphic to a closed interval
(Mostert \cite{Mostert57}, Bergery \cite{Bergery82}).
Orbits of interior points are principal and those of end points are singular.
Moreover principal orbits are hypersurfaces in $\tilde M$,
because $\tilde M/G$ is homeomorphic to a closed interval.
Suppose that $G(x)$ is a principal orbit which is a weakly reflective submanifold.
Then, by the slice representation theorem,
there exists a unit normal vector field $\xi$ on $G(x)$,
which is invariant under the action of $G$.
We take a geodesic $\gamma(t)$ of $\tilde M$
which satisfies an initial condition
$$
\gamma(0) = x, \qquad \gamma'(0) = \xi_x.
$$
Then $\gamma(t)$ is a section of the action of $G$ on $\tilde M$.
Since $dg_x(\xi_x) = \xi_{gx}$ for any $g \in G$,
$g\gamma(t)$ is a geodesic of $\tilde M$ which satisfies an initial condition
$$
g\gamma(0) = gx, \qquad (g\gamma)'(0) = \xi_{gx}.
$$
Since $G(x)$ is a weakly reflective submanifold of $\tilde M$,
there exists an isometry $\sigma$ of $\tilde M$ which satisfies
$$
\sigma(x) = x, \quad
d\sigma_x(\xi_x) = - \xi_x, \quad
\sigma(G(x)) = G(x),
$$
that is a reflection of $G(x)$ with respect to $\xi_x$.
We set
$$
G(x)_{\pm} = \{y \in G(x) \mid d\sigma_y(\xi_y) = \pm \xi_{\sigma(y)}\}.
$$
The sets $G(x)_+$ and $G(x)_-$ are closed subsets of $G(x)$,
and $G(x)$ is a disjoint union of $G(x)_+$ and $G(x)_-$
because $G(x)$ is a hypersurface in $\tilde M$.
Since $G(x)$ is connected and $x \in G(x)_-$,
we have $G(x) = G(x)_-$.
This implies that
$d\sigma_y(\xi_y) = - \xi_{\sigma(y)}$ for any $y \in G(x)$.
For any $g \in G$, $\sigma g\gamma(t)$ is a geodesic
which satisfies an initial condition
$$
\sigma g\gamma(0) = \sigma(gx), \qquad
(\sigma g\gamma)'(0) = d\sigma_{gx} (g\gamma)'(0)
= d\sigma_{gx}(\xi_{gx})
= - \xi_{\sigma gx}.
$$
Now we take $g_1 \in G$ such that $g_1x = \sigma(gx)$.
Then $\sigma g \gamma(t)$ and $g_1 \gamma(-t)$ are geodesics
of same initial conditions,
hence $\sigma g \gamma(t) = g_1 \gamma(-t) \in G(\gamma(-t))$.
Therefore we have $\sigma(G(\gamma(t))) \subset G(\gamma(-t))$ for each $t$.
Since $\sigma^{-1}$ is also a reflection of $G(x)$ at $x$,
we also have $\sigma^{-1}(G(\gamma(-t))) \subset G(\gamma(t))$
by the same discussion for $\sigma^{-1}$ and $\gamma(-t)$.
Thus $\sigma(G(\gamma(t))) = G(\gamma(-t))$.
This implies that two singular orbits can be expressed
as $G(\gamma(t_1))$ and $G(\gamma(-t_1))$ for some $t_1$.
Consequently we have the conclusion.
\end{proof}

\section{Orbits of $s$-representations}

A linear isotropy representation of a Riemannian symmetric pair
is called an $s$-representation as we mentioned in Introduction.
In the following sections, we will study orbits of $s$-representations
which are austere submanifolds and weakly reflective submanifolds.
For this purpose, we shall provide some fundamental notions
of orbits of $s$-representations in this section.

Let $G$ be a compact, connected Lie group
and $K$ a closed subgroup of $G$.
Assume that $\theta$ is an involutive automorphism of $G$
and $G_\theta^0 \subset K \subset G_\theta$, where
$$
G_\theta = \{g \in G \mid \theta(g) = g\}
$$
and $G_\theta^0$ is the identity component of $G_\theta$.
Then $(G, K)$ is a symmetric pair with respect to $\theta$.
We denote the Lie algebras of $G$ and $K$ 
by $\mathfrak g$ and $\mathfrak k$, respectively.
The involutive automorphism of $\mathfrak g$ induced from $\theta$
will be also denoted by $\theta$.
Then we have
$$
\mathfrak k = \{X \in \mathfrak g \mid \theta(X) = X\}.
$$
Take an inner product $\langle\; ,\; \rangle$ on $\mathfrak g$
which is invariant under $\theta$ and the adjoint representation of $G$.
Set
$$
\mathfrak m = \{X \in \mathfrak g \mid \theta(X) = - X\},
$$
then we have a canonical orthogonal direct sum decomposition
$$
\mathfrak g = \mathfrak k + \mathfrak m.
$$
Henceforth we assume that the symmetric pair $(G, K)$ is irreducible,
namely $K$ acts irreducibly on $\mathfrak m$.

Fix a maximal Abelian subspace $\mathfrak a$ in $\mathfrak m$
and a maximal Abelian subalgebra $\mathfrak t$ in $\mathfrak g$
containing $\mathfrak a$.
For $\alpha \in \mathfrak t$ we set
$$
\tilde{\mathfrak g}_\alpha
= \{X \in \mathfrak g^{\mathbf C} \mid
[H, X] = \sqrt{-1}\langle\alpha, H\rangle X\; (H \in \mathfrak t)\}
$$
and define the root system $\tilde R$ of $\mathfrak g$ by
$$
\tilde R = 
\{\alpha \in \mathfrak t - \{0\} \mid
\tilde{\mathfrak g}_\alpha \ne \{0\}\}.
$$
For $\alpha \in \mathfrak a$ we set
$$
\mathfrak g_\alpha
= \{X \in \mathfrak g^{\mathbf C} \mid
[H, X] = \sqrt{-1}\langle\alpha, H\rangle X\; (H \in \mathfrak a)\}
$$
and define the restricted root system $R$ of $(\mathfrak g, \mathfrak k)$ by
$$
R = \{\alpha \in \mathfrak a - \{0\} \mid
\mathfrak g_\alpha \ne \{0\}\}.
$$
Set
$$
\tilde R_0 = \tilde R \cap \mathfrak k
$$
and denote the orthogonal projection from $\mathfrak t$ to $\mathfrak a$ by
$H \mapsto \bar H$.
Then we have
$$
R = \{\bar\alpha \mid \alpha \in \tilde R - \tilde R_0\}.
$$
We take a basis of $\mathfrak t$ extended from a basis of $\mathfrak a$
and define the lexicographic orderings $>$ on $\mathfrak a$ and $\mathfrak t$
with respect to these bases.
Then for $H \in \mathfrak t$, $\bar H > 0$ implies $H > 0$.
We denote by $\tilde F$ the fundamental system of $\tilde R$
with respect to the ordering $>$.
Set
$$
\tilde F_0 = \tilde F \cap \tilde R_0,
$$
then the fundamental system $F$ of $R$ with respect to the ordering $>$
is given by
$$
F = \{\bar\alpha \mid \alpha \in \tilde F - \tilde F_0\}.
$$
We set
$$
\tilde R_+ = \{\alpha \in \tilde R \mid \alpha > 0\}, \qquad
R_+ = \{\alpha \in R \mid \alpha > 0\}.
$$
Then we have
$$
R_+ = \{\bar\alpha \mid \alpha \in \tilde R_+ - \tilde R_0\}.
$$
We also set
$$
\mathfrak k_0 = \{X \in \mathfrak k \mid [X, H] = 0\; (H \in \mathfrak a)\},
$$
and define
$$
\mathfrak k_\alpha
= \mathfrak k \cap (\mathfrak g_\alpha + \mathfrak g_{-\alpha}), \qquad
\mathfrak m_\alpha
= \mathfrak m \cap (\mathfrak g_\alpha + \mathfrak g_{-\alpha})
$$
for $\alpha \in R_+$.
Under these notations, we have the following lemma.

\begin{lem}[\cite{Takeuchi94}] \label{lem:basis}
\begin{enumerate}
\item We have orthogonal direct sum decompositions
$$
\mathfrak k = \mathfrak k_0 + \sum_{\alpha \in R_+} \mathfrak k_\alpha,
\qquad
\mathfrak m = \mathfrak a + \sum_{\alpha \in R_+} \mathfrak m_\alpha.
$$
\item
For each $\alpha \in \tilde R_+ - \tilde R_0$,
there exist $S_\alpha \in \mathfrak k$ and $T_\alpha \in \mathfrak m$
such that
$$
\{S_\alpha \mid \alpha \in \tilde R_+,\; \bar\alpha = \lambda\},
\quad
\{T_\alpha \mid \alpha \in \tilde R_+,\; \bar\alpha = \lambda\}
$$
are respectively orthonormal bases of $\mathfrak k_\lambda$
and $\mathfrak m_\lambda$
and that for $H \in \mathfrak a$
\begin{eqnarray*}
&&
[H, S_\alpha] = \langle\alpha, H\rangle T_\alpha, \quad
[H, T_\alpha] = - \langle\alpha, H\rangle S_\alpha, \quad
[S_\alpha, T_\alpha] = \bar\alpha,\\
&&
\mathrm{Ad}(\exp H)S_\alpha
= \cos\langle\alpha,H\rangle S_\alpha
+ \sin\langle\alpha,H\rangle T_\alpha,\\
&&
\mathrm{Ad}(\exp H)T_\alpha
= -\sin\langle\alpha,H\rangle S_\alpha
+\cos\langle\alpha,H\rangle T_\alpha.
\end{eqnarray*}
\end{enumerate}
\end{lem}

We define a subset $D$ of $\mathfrak a$ by
$$
D = \bigcup_{\alpha \in R}
\{H \in \mathfrak a \mid \langle\alpha, H\rangle = 0\}.
$$
A connected component of $\mathfrak a - D$ is a Weyl chamber.
We set
$$
C = \{H \in \mathfrak a \mid
\langle\alpha, H\rangle > 0\; (\alpha \in F)\}.
$$
Then $C$ is an open convex subset of $\mathfrak a$ and
the closure of $C$ is given by
$$
\bar C = \{H \in \mathfrak a \mid
\langle\alpha, H\rangle \ge 0\; (\alpha \in F)\}.
$$
For a subset $\Delta \subset F$,
we define
$$
C^\Delta = 
\{H \in \bar C \mid \langle\alpha, H\rangle > 0\; (\alpha \in \Delta),\;
\langle\beta, H\rangle = 0\; (\beta \in F-\Delta)\}.
$$

\begin{lem} \label{lem:docomposition of Weyl}
\begin{enumerate}
\item
For $\Delta_1\subset F$, the decomposition
$$
\overline{C^{\Delta_1}}=\bigcup_{\Delta\subset\Delta_1}C^\Delta
$$
is a disjoint union.
In particular,
$\displaystyle\bar C = \bigcup_{\Delta \subset F} C^\Delta$
is a disjoint union.
\item
For $\Delta_1, \Delta_2 \subset F$, $\Delta_1 \subset \Delta_2$
if and only if $C^{\Delta_1} \subset \overline{C^{\Delta_2}}$.
\end{enumerate}
\end{lem}

For each $\alpha \in F$,
we take $H_\alpha \in \mathfrak a$ such that
$$
\langle H_\alpha, \beta \rangle =
\left\{
\begin{array}{ll}
1 & (\beta = \alpha), \\
0 & (\beta \neq \alpha)
\end{array}
\right.
\quad (\beta \in F).
$$
Then we have
$$
\bar C =
\left\{\left.
\sum_{\alpha \in F}t_\alpha H_\alpha\; 
\right|\;
t_\alpha \ge 0
\right\},
$$
and for $\Delta \subset F$
$$
C^\Delta =
\left\{\left.
\sum_{\alpha \in \Delta}t_\alpha H_\alpha\; 
\right|\;
t_\alpha > 0
\right\}.
$$

We set
\begin{eqnarray*}
R^\Delta
& = &
R \cap (F - \Delta)_{\mathbf Z}, \\
R^\Delta_+
& = &
R^\Delta \cap R_+, \\
\mathfrak g^\Delta
& = &
\mathfrak k_0 + \mathfrak a 
+ \sum_{\alpha \in R^\Delta_+} (\mathfrak k_\alpha + \mathfrak m_\alpha).
\end{eqnarray*}
We also set
\begin{eqnarray*}
\mathfrak k^\Delta
& = &
\mathfrak g^\Delta \cap \mathfrak k
= \mathfrak k_0 + \sum_{\alpha \in R^\Delta_+} \mathfrak k_\alpha, \\
\mathfrak m^\Delta
& = &
\mathfrak g^\Delta \cap \mathfrak m
= \mathfrak a + \sum_{\alpha \in R^\Delta_+} \mathfrak m_\alpha.
\end{eqnarray*}
Then we have an orthogonal direct sum decomposition
$$
\mathfrak g^\Delta = \mathfrak k^\Delta + \mathfrak m^\Delta.
$$
For $H \in \mathfrak m$ we set
$$
Z^H_K = \{k \in K \mid \mathrm{Ad}(k)H = H\}.
$$
Then $Z^H_K$ is a closed subgroup of $K$ and the orbit $\mathrm{Ad}(K)H$ is
diffeomorphic to the coset manifold $K/Z_K^H$. 

Under these notations, we have the following lemma.

\begin{lem}[\cite{HKT00}] \label{lem:determination of R_+(H) from H}
Fix a subset $\Delta \subset F$.
For $H \in C^\Delta$ we have the following:
\begin{enumerate}
\item
$R^\Delta_+ = \{\alpha \in R_+\; |\; \langle\alpha, H\rangle = 0\}$,
\item
$R^\Delta = \{\alpha \in R\; |\; \langle\alpha, H\rangle = 0\}$,
\item
$\mathfrak g^\Delta = \{X \in \mathfrak g\; |\; [H, X] = 0\}$,
\item
$(\mathfrak g^\Delta, \mathfrak k^\Delta)$ is a symmetric pair
and its canonical decomposition is given by
$$
\mathfrak g^\Delta = \mathfrak k^\Delta + \mathfrak m^\Delta,
$$
\item $\mathfrak k^\Delta$ is the Lie algebra of $Z^H_K$.
\end{enumerate}
\end{lem}

Now we shall study an orbit $\mathrm{Ad}(K)H$
of the linear isotropy representation of $(G,K)$ through $H\in\mathfrak m$.
An orbit  $\mathrm{Ad}(K)H$ is a  submanifold of the hypersphere $S$
of radius $\| H \|$ in $\mathfrak m$.
From \cite{HKT00}, $\mathrm{Ad}(K)H$ is connected.
Since
$$
\mathfrak m=\bigcup_{k\in K}\mathrm{Ad}(k)\bar C,
$$
without loss of generalities we may assume $H\in\bar C$.
Moreover, from Lemma \ref{lem:docomposition of Weyl},
there exists $\Delta \subset F$ such that $H \in C^\Delta$.
For $X \in \mathfrak k$ we define a vector field $X^*$ on $\mathfrak m$ by
$$
X^*_x = \left.\frac{d}{dt}\right|_{t=0} \mathrm{Ad}(\exp tX)x = [X,x]
$$
at $x \in \mathfrak m$.
Then $X^*|_{\mathrm{Ad}(K)H}$ is a tangent vector field on $\mathrm{Ad}(K)H$.
From Lemma~\ref{lem:basis} we have the following lemma.

\begin{lem} \label{lem:orbits of s-representation}
For $\Delta \subset F$ and $H \in C^\Delta$,
the tangent space $T_H(\mathrm{Ad}(K)H)$
of the orbit $\mathrm{Ad}(K)H$ at $H$ and
the normal space $T_H^\perp(\mathrm{Ad}(K)H)$ in the hypersphere
can be expressed as
\begin{eqnarray*}
T_H(\mathrm{Ad}(K)H)
&=& \sum_{\alpha \in R_+ - R_+^\Delta} \mathfrak m_\alpha, \\
T_H^\perp(\mathrm{Ad}(K)H)
&=& H^\perp \cap \mathfrak a + \sum_{\alpha \in R_+^\Delta} \mathfrak m_\alpha
= \bigcup_{k \in Z^H_K} \mathrm{Ad}(k) (H^\perp \cap \mathfrak a).
\end{eqnarray*}
Let $h$ denote the second fundamental form of $\mathrm{Ad}(K)H$ at $H$
in the hypersphere $S$.
Then we have
$$
h(X_H^*,Y_H^*)
= [Y,[X,H]]^N,
$$
where $[Y,[X,H]]^N$ is $T_H^\perp(\mathrm{Ad}(K)H)$-component
of $[Y,[X,H]]$.
\end{lem}

\begin{proof}
From Lemma~\ref{lem:basis} we have
\begin{eqnarray*}
T_H(\mathrm{Ad}(K)H)
&=& \sum_{\alpha \in R_+ - R_+^\Delta} \mathfrak m_\alpha, \\
T_H^\perp(\mathrm{Ad}(K)H)
&=& H^\perp \cap \mathfrak a + \sum_{\alpha \in R_+^\Delta} \mathfrak m_\alpha
= H^\perp \cap \mathfrak m^\Delta
\end{eqnarray*}
Moreover, from Lemma~\ref{lem:determination of R_+(H) from H}
$$
\mathfrak m^\Delta = \bigcup_{k \in Z^H_K} \mathrm{Ad}(k) \mathfrak a.
$$
Since $\mathrm{Ad}(k)H = H$ for $k \in Z^H_K$, we have
$$
H^\perp \cap \mathfrak m^\Delta
= \bigcup_{k \in Z^H_K} \mathrm{Ad}(k) (H^\perp \cap \mathfrak a).
$$

The calculation of $X^*_x$ mentioned above
shows the representation of the second fundamental form.
\end{proof}

For orbits of $s$-representations which are minimal submanifolds
in the hypersphere,
the following theorem is known.

\begin{thm}[\cite{HSTT00}] \label{thm:minimal orbit}
Fix a hypersphere $S$ in $\mathfrak m$ centered at $0$.
For each subset $\Delta\subset F$,
there exists a unique $H\in C^\Delta \cap S$ such that
$\mathrm{Ad}(K)H$ is a minimal submanifold in $S$.
\end{thm}

However, in general we can not determine $H$
where $\mathrm{Ad}(K)H$ is a minimal submanifold in $S$ explicitly.
In the following two sections,
we will give the complete lists of $H$
where $\mathrm{Ad}(K)H$ is an austere submanifold
and a weakly submanifold in $S$.

\section{Weakly reflective orbits of $s$-representations}
\label{sec:wr-orbit}

In this section, we shall study orbits of irreducible $s$-representations
which are weakly reflective submanifolds in the hypersphere.
In the next section, we will study austere orbits.
Since these two properties of orbits are invariant
under scalar multiples on the vector spaces,
we do not discriminate the difference of the length of a vector.
The following theorem is the main result of this section.
We shall follow the notations of root systems in \cite{Bourbaki75}.

\begin{thm} \label{thm:weakly reflective orbit}
An orbit of an irreducible $s$-representation which is a weakly reflective
submanifold in the hypersphere is one of the following list:
\begin{enumerate}
\item
an orbit through a restricted root vector
{\rm (Proposition~\ref{pro:orbit of root})},
\item
the orbit through the vector $2e_1 - e_2 - e_3$ or $e_1 + e_2 - 2e_3$
of the linear isotropy representation of a compact symmetric pair
with restricted root system $\{\pm(e_i-e_j)\}$ of type $A_2$
{\rm (Proposition~\ref{pro:orbit of A_2})},
\item
the orbit through the vector $e_1 + e_2 - e_3 - e_4$
of the linear isotropy representation of a compact symmetric pair
with restricted root system $\{\pm(e_i-e_j)\}$ of type $A_3$
{\rm (Proposition~\ref{pro:orbit of A_3})},
\item
the orbit through the vector $e_1$
of the linear isotropy representation of a compact symmetric pair
with restricted root system $\{\pm(e_i \pm e_j)\}$ of type $D_l$
{\rm (Proposition~\ref{pro:e_1-orbit of D_l})},
\item
the orbit through the vector $e_1 + e_2 + e_3 \pm e_4$
of the linear isotropy representation of a compact symmetric pair
with restricted root system $\{\pm(e_i \pm e_j)\}$ of type $D_4$
{\rm (Proposition~\ref{pro:orbit of D_4})}.
\end{enumerate}
\end{thm}
Here we prove that orbits listed above are
weakly reflective submanifolds.
In the next section, we will classify all austere orbits
of irreducible $s$-representations and show that all weakly reflective orbits
can be obtained in Theorem~\ref{thm:weakly reflective orbit}.
For this purpose, we shall first give the following lemma.

\begin{lem} \label{lem:criterion for weakly reflective orbits}
For $H \in \mathfrak a$, the orbit $\mathrm{Ad}(K)H$
is a weakly reflective submanifold in the hypersphere $S$ if and only if,
for any $\xi \in H^\perp \cap \mathfrak a$,
there exists a linear isometry $\sigma_\xi$ of $\mathfrak m$
which satisfies
\begin{equation} \label{eq:criterion for weakly reflective orbits}
\sigma_\xi(H) = H, \qquad
\sigma_\xi(\xi) = - \xi, \qquad
\sigma_\xi(\mathrm{Ad}(K)H) = \mathrm{Ad}(K)H.
\end{equation}
\end{lem}

\begin{proof}
From Lemma \ref{lem:orbits of s-representation},
the normal space of the orbit $\mathrm{Ad}(K)H$ at $H$ in $S$ is given by
$$
T_H^\perp (\mathrm{Ad}(K)H)
= H^\perp \cap \mathfrak a
+ \sum_{\alpha \in R_+^\Delta} \mathfrak m_\alpha
= \bigcup_{k \in Z^H_K} \mathrm{Ad}(k) (H^\perp \cap \mathfrak a).
$$
If $\mathrm{Ad}(K)H$ is a weakly reflective submanifold in $S$,
then for $\xi \in H^\perp \cap \mathfrak a$
there exists a linear isometry $\sigma_\xi$ of $\mathfrak m$ which satisfies
$$
\sigma_\xi(H) = H, \qquad
(d\sigma_\xi)_H \xi = - \xi, \qquad
\sigma_\xi(\mathrm{Ad}(K)H) = \mathrm{Ad}(K)H.
$$
Here we have $(d\sigma_\xi)_H = \sigma_\xi$,
since $\sigma_\xi$ is a linear isometry.

Conversely, assume that $\mathrm{Ad}(K)H$ satisfies
the condition (\ref{eq:criterion for weakly reflective orbits}).
We take an arbitrary normal vector $\xi \in T_H^\perp (\mathrm{Ad}(K)H)$.
From Lemma~\ref{lem:orbits of s-representation},
there exists $k_0 \in Z_K^H$ such that
$\mathrm{Ad}(k_0)\xi \in H^\perp \cap \mathfrak a$.
Then, from the assumption,
there exists a linear isometry $\sigma$ which satisfies
$$
\sigma(H) = H, \qquad
\sigma \mathrm{Ad}(k_0)\xi = - \mathrm{Ad}(k_0)\xi, \qquad
\sigma(\mathrm{Ad}(K)H) = \mathrm{Ad}(K)H.
$$
We now define
$\sigma_\xi = \mathrm{Ad}(k_0)^{-1}\sigma\mathrm{Ad}(k_0)$.
Then $\sigma_\xi$ satisfies
$$
\sigma_\xi(H) = H, \qquad
\sigma_\xi(\xi) = - \xi, \qquad
\sigma_\xi(\mathrm{Ad}(K)H) = \mathrm{Ad}(K)H.
$$
Thus $\sigma_\xi$ is a reflection of $\mathrm{Ad}(K)H$
with respect to a normal vector $\xi$ at $H$.
Since $\mathrm{Ad}(K)H$ is a homogeneous submanifold,
we have a reflection with respect to any normal vector at arbitrary point.
Consequently $\mathrm{Ad}(K)H$ is a weakly reflective submanifold in $S$.
\end{proof}

\begin{pro}\label{pro:orbit of root}
An orbit through a restricted root vector of the linear isotropy representation
of an irreducible compact symmetric pair
is a weakly reflective submanifold in the hypersphere $S$.
\end{pro}

\begin{proof}
Let $\alpha_0$ be a restricted root vector and put $H = \alpha_0$.
The reflection $s_{\alpha_0}$ on $\mathfrak a$ with respect to $\alpha_0$
is given by
$$
s_{\alpha_0}(X)
= X - \frac{2\langle\alpha_0, X\rangle}
{\langle\alpha_0, \alpha_0\rangle} \alpha_0
\qquad (X \in \mathfrak a)
$$
and satisfies
$$
s_{\alpha_0}(H) = -H, \qquad
s_{\alpha_0}|_{\mathfrak a \cap H^\perp} = 1_{\mathfrak a \cap H^\perp}.
$$
The reflection $s_{\alpha_0}$ is an element of the Weyl group,
hence there exists $k_0 \in N_K$ such that
$\mathrm{Ad}(k_0)|_{\mathfrak a} = s_{\alpha_0}$, where
$$
N_K = \{k \in K \mid \mathrm{Ad}(k)\mathfrak a = \mathfrak a\}.
$$
Therefore
$$
-H = \mathrm{Ad}(k_0)H \in \mathrm{Ad}(K)H,
$$
and we have $\mathrm{Ad}(K)(-H) = \mathrm{Ad}(K)H$.
We define a linear isometry $\sigma$ of $\mathfrak m$ by
$$
\sigma = -\mathrm{Ad}(k_0)|_{\mathfrak m}.
$$
Then $\sigma$ satisfies
$$
\sigma(H) = H, \qquad
\sigma|_{\mathfrak a \cap H^\perp} = -1|_{\mathfrak a \cap H^\perp}, \qquad
\sigma(\mathrm{Ad}(K)H) = \mathrm{Ad}(K)H.
$$
Thus, from Lemma~\ref{lem:criterion for weakly reflective orbits},
$\mathrm{Ad}(K)H$ is a weakly reflective submanifold in $S$.
\end{proof}

\begin{pro}\label{pro:orbit of A_2}
The orbit through the vector $2e_1 - e_2 - e_3$ or $e_1 + e_2 - 2e_3$
of the linear isotropy representation of a compact symmetric pair
with restricted root system $\{\pm(e_i - e_j)\}$ of type $A_2$
is a weakly reflective submanifold in the hypersphere $S$.
\end{pro}

\begin{proof}
Since the symmetric pair $(G,K)$ is of rank $2$,
the action of $K$ on $S$ is cohomogeneity one.
The vector $2e_1 - e_2 - e_3$ (resp. $e_1 + e_2 - 2e_3$) is
orthogonal to a restricted root $e_2 - e_3$ (resp. $e_1 - e_2$).
Therefore the orbit of $K$
through $2e_1 - e_2 - e_3$ (resp. $e_1 + e_2 - 2e_3$) is a singular orbit.
Hence from Proposition~\ref{pro:cohomogeneity one action},
this orbit is a weakly reflective submanifold in $S$.
\end{proof}

\begin{pro} \label{pro:orbit of A_3}
The orbit through the vector $e_1 + e_2 - e_3 - e_4$
of the linear isotropy representation of a compact symmetric pair
with restricted root system $\{\pm(e_i - e_j)\}$ of type $A_3$
is a weakly reflective submanifold in the hypersphere $S$.
\end{pro}

\begin{proof}
Put $H = e_1 + e_2 - e_3 - e_4$.
The set $R^\Delta_+$ of all positive restricted roots which are orthogonal
to $H$ is given by
$$
R^\Delta_+ = \{ e_1 - e_2,\ e_3 - e_4 \}.
$$
Let $s_{e_1 - e_2}$ and $s_{e_3 - e_4}$ be the reflections
with respect to restricted roots $e_1 - e_2$ and $e_3 - e_4$, respectively.
Then $s_{e_1 - e_2}$ and $s_{e_3 - e_4}$ are elements of the Weyl group,
hence there exist $k_0, k_1 \in N_K$ such that
$$
s_{e_1 - e_2} = \mathrm{Ad}(k_0)|_{\mathfrak a}, \qquad
s_{e_3 - e_4} = \mathrm{Ad}(k_1)|_{\mathfrak a}.
$$
We now define a linear isometry of $\mathfrak m$ by
$$
\sigma(X) = \mathrm{Ad}(k_0) \mathrm{Ad}(k_1) X
\qquad (X \in \mathfrak m).
$$
Then $\sigma$ satisfies
$$
\sigma(H) = H, \qquad
\sigma|_{\mathfrak a \cap H^\perp} = -1_{\mathfrak a \cap H^\perp}, \qquad
\sigma(\mathrm{Ad}(K)H) = \mathrm{Ad}(K)H.
$$
Thus from Lemma~\ref{lem:criterion for weakly reflective orbits},
$\mathrm{Ad}(K)H$ is a weakly reflective submanifold in $S$.
\end{proof}

\begin{pro}\label{pro:e_1-orbit of D_l}
The orbit through the vector $e_1$
of the linear isotropy representation of a compact symmetric pair
with restricted root system $\{\pm e_i \pm e_j\}$ of type $D_l$
is a weakly reflective submanifold in the hypersphere $S$.
\end{pro}

\begin{proof}
An irreducible compact symmetric pair with restricted root system of
type $D_l$ is one of $(SO(2l)\times SO(2l), SO(2l)^*)$ and
$(SO(2l), SO(l)\times SO(l))$.

First we consider the case of $(SO(2l)\times SO(2l), SO(2l)^*)$.
In this case, $\mathfrak m$ can be identified with $\mathfrak o(2l)$
in a natural manner.
We take a maximal Abelian subalgebra
$$
\mathfrak a =
\left\{\left.
\mathrm{diag}\left\{
\left[\begin{array}{cc}
0 & -t_1 \\ t_1 & 0
\end{array}\right],
\dots,
\left[\begin{array}{cc}
0 & -t_l \\ t_l & 0
\end{array}\right]
\right\}\;
\right|\;
t_1, \dots, t_l \in \mathbf R
\right\}
$$
of $\mathfrak o(2l)$, and put
$$
H = e_1 =
\mathrm{diag}\left\{
\left[\begin{array}{cc}
0 & -1 \\ 1 & 0
\end{array}\right],
0, \dots, 0
\right\}.
$$
We define a linear isometry $\sigma$ of $\mathfrak o(2l)$ by
$$
\sigma(X) = sXs \qquad (X \in \mathfrak o(2l)),
$$
where
$$
s =
\mathrm{diag}\left\{
\left[\begin{array}{cc}
1 & 0 \\ 0 & 1
\end{array}\right],
\left[\begin{array}{cc}
-1 & 0 \\ 0 & 1
\end{array}\right],
\dots,
\left[\begin{array}{cc}
-1 & 0 \\ 0 & 1
\end{array}\right]
\right\}
\in O(2l).
$$
Then $\sigma$ is an isometry of $S$ and satisfies
$$
\sigma(H) = H, \qquad
\sigma|_{\mathfrak a \cap H^\perp} = -\mathrm{id}_{\mathfrak a \cap H^\perp},
\qquad
\sigma(\mathrm{Ad}(K)H) = \mathrm{Ad}(K)H.
$$
Hence from Lemma~\ref{lem:criterion for weakly reflective orbits},
$\mathrm{Ad}(K)H$ is a weakly reflective submanifold in $S$.

Second, we consider the case of $(SO(2l), SO(l)\times SO(l))$.
We take a maximal Abelian subspace
$$
\mathfrak a
= \left\{
\left.
\left[\begin{array}{cc}
0 & X \\ - X & 0
\end{array}\right]
\;\right|\;
X = \mathrm{diag}(t_1, \dots, t_l),\; t_i \in \mathbf R
\right\},
$$
and put
$$
H = e_1
= \left[\begin{array}{cc}
0 & X_0 \\ -X_0 & 0
\end{array}\right]
\in \mathfrak o(2l),
$$
where
$$
X_0 =
\mathrm{diag}\{1, 0, \dots, 0\}
\in M_l(\mathbf R).
$$
We define a linear isometry $\sigma$ of $\mathfrak m$ by
$$
\sigma(X) =
\left[\begin{array}{cc}
s & 0 \\ 0 & I_l
\end{array}\right]
X
\left[\begin{array}{cc}
s & 0 \\ 0 & I_l
\end{array}\right]
\qquad (X \in \mathfrak m),
$$
where
$$
s =
\mathrm{diag}\{1, -1, \dots, -1\}
\in O(l).
$$
Then $\sigma$ is an isometry of $S$ and satisfies
$$
\sigma(H) = H, \qquad
\sigma|_{\mathfrak a \cap H^\perp} = -\mathrm{id}_{\mathfrak a \cap H^\perp},
\qquad
\sigma(\mathrm{Ad}(K)H) = \mathrm{Ad}(K)H.
$$
Hence from Lemma~\ref{lem:criterion for weakly reflective orbits},
$\mathrm{Ad}(K)H$ is a weakly reflective submanifold in $S$.
\end{proof}

\begin{pro} \label{pro:orbit of D_4}
The orbit through the vector $e_1 + e_2 + e_3 \pm e_4$
of the linear isotropy representation of a compact symmetric pair
with restricted root system $\{\pm e_i \pm e_j\}$ of type $D_4$
is a weakly reflective submanifold in the hypersphere $S$.
\end{pro}

\begin{proof}
We take a fundamental system of $\{\pm e_i \pm e_j\}$ of type $D_4$:
$$
\alpha_1 = e_1 - e_2, \quad
\alpha_2 = e_2 - e_3, \quad
\alpha_3 = e_3 - e_4, \quad
\alpha_4 = e_3 + e_4.
$$
The automorphism group of the Dynkin diagram is
the permutation group of $\{\alpha_1, \alpha_3, \alpha_4\}$.
So there exists an automorphism of $\mathfrak m$
mapping $\alpha_1$ to $\alpha_4$ and fixing $\alpha_3$,
which gives an equivalence of the orbits through $e_1$
and $e_1 + e_2 + e_3 + e_4$.
Thus the orbit through $e_1 + e_2 + e_3 + e_4$ is also a weakly reflective
submanifold in the hypersphere.
Similarly the permutaion of $\alpha_1$ and $\alpha_3$
gives an equivalence of the orbits through $e_1$
and $e_1 + e_2 + e_3 - e_4$.
Thus the orbit through $e_1 + e_2 + e_3 - e_4$ is also a weakly reflective
submanifold in the hypersphere.
\end{proof}

In the case of $(SO(8), SO(4)\times SO(4))$
we can explicitly represent a reflection of the orbit
though $e_1 + e_2 + e_3 + e_4$.
The linear isotropy representation is equivalent to
$$
(g_1, g_2)\cdot X = g_1Xg_2^{-1}
\qquad ((g_1, g_2) \in SO(4)\times SO(4),\; X \in M_4(\mathbf R)).
$$
Let $e_i$ denotes an element of $M_4(\mathbf R)$
whose $(i,i)$ component is $1$ and others are $0$.
Then the orbit through $e_1 + e_2 + e_3 + e_4$ is $SO(4)$
in $M_4(\mathbf R)$.

For $z_1\otimes z_2 \in \mathbf H\otimes\mathbf H$,
we define
$$
\phi_{z_1\otimes z_2} : \mathbf H \to \mathbf H
\; ;\;
z \mapsto z_1z\bar z_2.
$$
Since $\mathbf H \cong \mathbf R^4$,
we can regard $\phi_{z_1\otimes z_2}$ as an element of $M_4(\mathbf R)$
and $\phi$ induces an isomorphism
$M_4(\mathbf R) \cong \mathbf H \otimes \mathbf H$ of real algebras.
We define an involutive isometry $\sigma$ of $\mathbf H\otimes\mathbf H$ by
$$
\sigma : \mathbf H\otimes\mathbf H \to \mathbf H\otimes\mathbf H
\; ;\;
z_1\otimes z_2 \mapsto z_1\otimes\bar z_2.
$$
We also denote by $\sigma$ the linear isometry of
$M_4(\mathbf R)$ induced from $\sigma$ through $\phi$.
We note that
$$
SO(4) =
\{\phi_{z_1\otimes z_2} \mid z_1, z_2 \in Sp(1)\}.
$$
Moreover
$$
\{z_1\otimes z_2 \mid z_1, z_2 \in Sp(1)\}
\subset \mathbf H\otimes\mathbf H
$$
is invariant under $\sigma$.
Therefore $SO(4)$ is invariant under $\sigma$.
The identity element $I$ is fixed by the action of $\sigma$.
The normal space of $SO(4)$ at $I$ in $S^{15}(2)$ is given by
$$
T_I^\perp(SO(4)) =
\{X \in M_4(\mathbf R) \mid X : \mbox{symmetric},\ \mathrm{tr}(X) = 0\}.
$$
It is easy to see that $T_I^\perp(SO(4))$ is contained in the eigenspace
of $\sigma$ for an eigenvalue $-1$.
Thus $\sigma$ is a reflection of $SO(4)$
with respect to an arbitrary normal vector at $I$.

\section{Austere orbits of $s$-representations}
\label{sec:austere-orbit}

In this section we classify all orbits of irreducible $s$-representations
which are austere submanifolds in the hypersphere $S$.
In the previous section we showed that all orbits
through a restricted root vector (or its scalar multiple)
are weakly reflective, hence austere.
Therefore, hereafter we shall concern with other orbits.
We will also determine austere orbits
which are not weakly reflective submanifolds.
Then we will complete to prove Theorem~\ref{thm:weakly reflective orbit}.

The classification of austere orbits is following:

\begin{thm} \label{thm:austere}
An orbit of an irreducible $s$-representation which is an austere
submanifold in the hypersphere is one of the following list:
\begin{enumerate}
\item
an orbit through a restricted root vector,
\item
the orbit through the vector $2e_1 - e_2 - e_3$ or $e_1 + e_2 - 2e_3$
of the linear isotropy representation of a compact symmetric pair
with restricted root system $\{\pm(e_i-e_j)\}$ of type $A_2$
{\rm (Proposition~\ref{pro:orbit of A_2} or \ref{pro:austere-A})},
\item
the orbit through the vector $e_1 + e_2 - e_3 - e_4$
of the linear isotropy representation of a compact symmetric pair
with restricted root system $\{\pm(e_i-e_j)\}$ of type $A_3$
{\rm (Proposition~\ref{pro:austere-A})},
\item
the orbit through the vector $e_1$
of the linear isotropy representation of a compact symmetric pair
with restricted root system $\{\pm(e_i \pm e_j)\}$ of type $D_l$
{\rm (Proposition~\ref{pro:austere-D})},
\item
the orbit through the vector $e_1 + e_2 + e_3 \pm e_4$
of the linear isotropy representation of a compact symmetric pair
with restricted root system $\{\pm(e_i \pm e_j)\}$ of type $D_4$
{\rm (Proposition~\ref{pro:austere-D})},
\item
the orbit through the vector $e_1+\frac{e_1+e_2}{\sqrt{2}}$
of the linear isotropy representation of a compact symmetric pair
with restricted root system $\{\pm e_i, \pm e_i \pm e_j\}$ of type $B_2$
whose multiplicities are constant
{\rm (Proposition~\ref{pro:austere-BC})},
\item
the orbit through the vector $\alpha_1+\frac{\alpha_2}{\sqrt{3}}$
of the linear isotropy representation of a compact symmetric pair
with restricted root system of type $G_2$
{\rm (Proposition~\ref{pro:austere-G})}.
\end{enumerate}
\end{thm}

\begin{rem}
In the case where the rank of the symmetric pair is equal to two,
any principal orbit of $s$-representations is
an isoparametric hypersurface in the hypersphere.
The family of isoparametric hypersurfaces has a unique
minimal isoparametric hypersurface.
The theorem above shows some of minimal isoparametric hypersurfaces
are austere, furthermore weakly reflective, and some of them are not.
\end{rem}

Before giving a proof of Theorem~\ref{thm:austere}
we shall provide some preliminaries.
Let $(G,K)$ be an irreducible compact symmetric pair.
We shall use the notations of previous sections.
From Lemma \ref{lem:orbits of s-representation},
for a normal vector
$\xi \in T_H^\perp(\mathrm{Ad}(K)H)$,
the shape operator $A_\xi$ of $\mathrm{Ad}(K)H$
in the hypersphere $S$ is given by
\begin{equation} \label{eq:5-2}
\langle A_\xi(X^*), Y^*\rangle
= \langle h(X^*, Y^*), \xi\rangle
= \langle[Y, [X, H]], \xi\rangle
= - \langle[X, H], [Y, \xi]\rangle.
\end{equation}
For simplification,
we discuss a normalization of a normal vector $\xi$.
From Lemma~\ref{lem:orbits of s-representation},
there exists $k \in Z_K^H$ such that
$\mathrm{Ad}(k)\xi \in H^\perp \cap \mathfrak a$.
Then
\begin{eqnarray*}
\langle A_\xi(X^*), Y^*\rangle
&=& \langle \mathrm{Ad}(k)h(X^*, Y^*), \mathrm{Ad}(k)\xi\rangle \\
&=& \langle h(\mathrm{Ad}(k)X^*, \mathrm{Ad}(k)Y^*),
\mathrm{Ad}(k)\xi \rangle \\
& = &
\langle A_{\mathrm{Ad}(k)\xi}\mathrm{Ad}(k)X^*, \mathrm{Ad}(k)Y^*
\rangle \\
& = &
\langle \mathrm{Ad}(k)^{-1}A_{\mathrm{Ad}(k)\xi}(\mathrm{Ad}(k)X^*), Y^*
\rangle.
\end{eqnarray*}
Thus we have
$$
A_\xi = \mathrm{Ad}(k)^{-1}A_{\mathrm{Ad}(k)\xi}\mathrm{Ad}(k).
$$
This implies that eigenvalues of $A_{\mathrm{Ad}(k)\xi}$
and their multiplicities coincide with those of $A_\xi$.
Hence, in order to show whether an orbit $\mathrm{Ad}(K)H$ is austere,
it suffices to check eigenvalues of $A_\xi$
for $\xi \in H^\perp \cap \mathfrak a$.
Hereafter we assume that $\xi \in H^\perp \cap \mathfrak a$.

From Lemmas \ref{lem:basis} and \ref{lem:orbits of s-representation}
$$
\{T_\alpha \mid \alpha \in \tilde R_+ - \tilde R_+^\Delta\}
$$
is an orthonormal basis of $T_H(\mathrm{Ad}(K)H)$.
For $\alpha, \beta \in \tilde R_+ - \tilde R_+^\Delta$ we have
$$
\langle A_\xi((S_\alpha^*)_H), (S_\beta^*)_H\rangle
= \langle\alpha, H\rangle\langle\beta, H\rangle
\langle A_\xi(T_\alpha), T_\beta\rangle.
$$
On the other hand, from (\ref{eq:5-2}), we have
$$
\langle A_\xi((S_\alpha^*)_H), (S_\beta^*)_H \rangle
= - \langle [S_\alpha, H], [S_\beta, \xi] \rangle
= - \langle \alpha, H \rangle \langle \beta, \xi \rangle \delta_{\alpha\beta}.
$$
Therefore we have
$$
A_\xi(T_\alpha)
= - \frac{\langle\alpha, \xi\rangle}{\langle\alpha, H\rangle}T_\alpha.
$$
This shows that $T_\alpha$ is an eigenvector of $A_\xi$
and its eigenvalue is
$$
- \frac{\langle\alpha, \xi\rangle}{\langle\alpha, H\rangle}.
$$
Hence $\mathrm{Ad}(K)H$ is an austere submanifold in $S$
if and only if, for any $\xi \in H^\perp \cap \mathfrak a$, the set
$$
\left\{\left.
\frac{\langle\alpha, \xi\rangle}{\langle\alpha, H\rangle}
\;\right|\;
\alpha \in \tilde R_+ - \tilde R_+^\Delta
\right\}
$$
is symmetric (concerning multiplicities) by the multiplication of $-1$.
We shall describe this condition
in terms of a finite subset of a Euclidean space.

Let $A$ be a finite subset of a Euclidean space $V$.
We consider a condition that, for any $v \in V$, 
$$
\{\langle a, v\rangle \mid a \in A\}
$$
is symmetric (concerning multiplicities) by the multiplication of $-1$.
This condition is equivalent to a condition that
$A$ is symmetric by the multiplication of $-1$ on $V$.
Indeed, it is obvious that 
$\{ \langle a,v \rangle \ | \ a \in A \}$ is symmetric whenever
$A$ is symmetric.
Conversely, fix an arbitrary $a \in A$.
From the assumption we have
$$
V= \bigcup_{b \in A}\{ v \in V\ | \ \langle a, v \rangle = -\langle b,
v\rangle \}.
$$
If $-a \notin A$, then the right hand side consists of finite union of
hyperplanes of $V$.
This is a contradiction.
Therefore $-a\in A$.
Consequently $A$ is symmetric by the multiplication of $-1$ on $V$.

Let $p_H : \mathfrak a \to H^\perp \cap \mathfrak a$ denote
the orthogonal projection.
An orbit $\mathrm{Ad}(K)H$ is austere in $S$
if and only if the set
$$
\left\{\left.
\frac{p_H(\alpha)}{\langle\alpha, H\rangle}
\;\right|\;
\alpha \in R_+ - R_+^\Delta
\right\}
$$
is symmetric (concerning multiplicities) by the multiplication of $-1$.
By this criterion,
we can easily see that orbits listed in Theorem~\ref{thm:austere}
are austere submanifolds in the hypersphere $S$.
Hereafter we shall prove that all austere orbits can be obtained
in Theorem~\ref{thm:austere}.

We set
$\mbox{\boldmath$R$}R=\{x\alpha\;|\;x\in \mbox{\boldmath$R$},\alpha\in R\}$.
We have alrady showed that the orbit through any element in
$\mbox{\boldmath$R$}R$ is weakly reflective in the hypersphere,
so we consider the orbits through elements in 
$\frak a -\mbox{\boldmath$R$}R$

\begin{lem} \label{lem:A}
For $H \in \frak a -\mbox{\boldmath$R$}R$,
the orbit $\mathrm{Ad}(K)H$ is an austere submanifold in $S$
if and only if 
there exist a mapping $f:R_+-R_+^\Delta\rightarrow R_+-R_+^\Delta$
without fixed points, and
constants $n_\alpha \neq 0, \epsilon_\alpha = \pm 1$
for each $\alpha\in R_+-R_+^\Delta$ such that
\begin{equation}\label{eqn:1}
H=n_\alpha\left(
\frac{\alpha}{\| \alpha \|}
+\epsilon_\alpha\frac{f(\alpha )}{\| f(\alpha ) \|}
\right),
\end{equation}
and
\begin{equation}
\label{eqn:multiplicity}
\sum_{\mu \in R_+ - R_+^\Delta \atop \mu \mathrel{/\!/} \alpha} m(\mu)
= \sum_{\nu \in R_+ - R_+^\Delta \atop \nu \mathrel{/\!/} f(\alpha)} m(\nu).
\end{equation}
Here we denote by $m(\mu)$ the multiplicity of a restricted root $\mu$.

Excepting the case where the restricted root system $R$ is of type $BC$,
the equality $(\ref{eqn:multiplicity})$ is equivalent to
$m(\alpha) = m(f(\alpha))$,
moreover $\# (R_+-R_+^\Delta )$ is even and $f^2=1$.
\end{lem}

\begin{proof}
The orthogonal projection $p_H$ is defined by
$$
p_H(X) = X - \frac{\langle X, H \rangle}{\langle H, H \rangle} H
\qquad (X \in \mathfrak a).
$$
Therefore $\mbox{Ad}(K)H$ is an austere submanifold in $S$
if and only if the set
$$
\left\{\left.
\frac{\alpha}{\langle\alpha ,H\rangle}-\frac{H}{\langle H,H\rangle}
\;\right|\; \alpha\in R_+-R_+^\Delta
\right\}
$$
is symmetric (concerning multiplicities) by the multiplication of $-1$.
In other words, there exists a mapping
$f:R_+-R_+^\Delta\rightarrow R_+-R_+^\Delta$ which satisfies
\begin{equation} \label{eq:5-4}
\frac{f(\alpha )}{\langle f(\alpha ),H\rangle}-\frac{H}{\langle H,H\rangle}=
-\frac{\alpha}{\langle\alpha ,H\rangle}+\frac{H}{\langle H,H\rangle}
\end{equation}
and
\begin{eqnarray*}
&&
\sum\left\{
m(\mu) \left|\;
\mu \in R_ + -R_+^\Delta,\
\frac{\mu}{\langle\mu, H\rangle} - \frac{H}{\langle H, H\rangle}
= \frac{\alpha}{\langle\alpha, H\rangle} - \frac{H}{\langle H, H\rangle}
\right.
\right\} \\
& = &
\sum\left\{
m(\nu) \left|\;
\nu \in R_ + -R_+^\Delta,\
\frac{\nu}{\langle\nu, H\rangle} - \frac{H}{\langle H, H\rangle}
= \frac{f(\alpha)}{\langle f(\alpha), H\rangle} 
- \frac{H}{\langle H, H\rangle}
\right.
\right\}
\end{eqnarray*}
for any $\alpha\in R_+-R_+^\Delta$.
This condition for the multiplicities is equivalent to (\ref{eqn:multiplicity}).
From (\ref{eq:5-4}), if $f$ has a fixed point $\alpha$,
then $H \in \mbox{\boldmath$R$}R$.
Thus $f$ has no fixed points.

If we assume (\ref{eq:5-4}), then there exist non-zero real numbers $x, y$
so that $H = x \alpha + y f(\alpha )$.
Applying this to the equation (\ref{eq:5-4}),
we have a quadratic equation
$$
\| f(\alpha ) \|^2 y^2 = \| \alpha \|^2 x^2
$$
with respect to $x$ and $y$.
Thus we have
$$
y = \pm\frac{\| \alpha \|}{\| f(\alpha ) \|} x,
$$
hence $H$ can be expressed as
$$
H = x \alpha \pm\frac{\|\alpha\|}{\|f(\alpha)\|}x f(\alpha)
= x\|\alpha\|\left(
\frac{\alpha}{\|\alpha\|} \pm \frac{f(\alpha)}{\|f(\alpha)\|}
\right).
$$
Since this equality holds for any $\alpha\in R_+-R_+^\Delta$,
we have the condition (\ref{eqn:1}).

Replacing $\alpha$ in (\ref{eq:5-4}) by $f(\alpha )$,
we have
\begin{equation} \label{eqn:3}
\frac{f^2(\alpha )}{\langle f^2(\alpha ),H\rangle}
+\frac{f(\alpha )}{\langle f(\alpha ),H\rangle}=\frac{2H}{\langle H,H\rangle}.
\end{equation}
From equations (\ref{eq:5-4}) and (\ref{eqn:3}), we have
$$
\frac{\alpha}{\langle \alpha ,H\rangle}
=\frac{f^2(\alpha )}{\langle f^2(\alpha ),H\rangle}.
$$
The above discussion stands for any restricted root systems $R$,
including of type $BC$.
Henceforth we assume that $R$ is not of type $BC$.
Then $\alpha$ is the only element of $R_+-R_+^\Delta$
which is a scalar multiple of $\alpha$.
Thus $f^2(\alpha ) = \alpha$.
Since $f$ has no fixed points,
$\# (R_+-R_+^\Delta )$ is even.
This completes the proof.
\end{proof}

\begin{pro} \label{pro:austere-A}
In the case where $R$ is of type $A_l$,
an austere orbit is one of the following
except orbits through a restricted root vector:
\begin{enumerate}
\item when $l=2$, the orbit through $H=2e_1-e_2-e_3$ or $e_1+e_2-2e_3$,
\item when $l=3$, the orbit through $H=e_1+e_2-e_3-e_4$.
\end{enumerate}
\end{pro}

\begin{proof}
In the case of $R=A_l$, 
$R_+$ is given by $R_+ = \{ e_i-e_j \mid i<j \}$.
Since all restricted roots have constant multiplicities,
the condition (\ref{eqn:multiplicity})
of Lemma~\ref{lem:A} is always satisfied.
From Lemma~\ref{lem:A},
without loss of generalities, we may assume
that $H=\pm\mbox{(a positive root)}\pm\mbox{(a positive root)}$.
Moreover
since any root can be translated to $e_1-e_2$ by the action of the Weyl group,
we may assume that $H=(e_1-e_2)\pm\mbox{(a positive root)}$.
The positive root in the second term of $H$ is one of
$e_1-e_i \; (3\leq i), e_2-e_j \; (3\leq i), e_i-e_j \; (3\leq i<j)$.

In the case of $H = (e_1-e_2)\pm (e_1-e_i) \; (3\leq i)$,
$e_i$ can be translated to $e_3$ by the action of an element of the Weyl group
which fixes both $e_1$ and $e_2$.
Therefore we can put
$$
H = (e_1-e_2)\pm (e_1-e_3)
= \left\{\begin{array}{l}
(e_1 - e_2) + (e_1 - e_3) = 2e_1 - e_2 - e_3 \\
(e_1 - e_2) - (e_1 - e_3) = - e_2 + e_3 \quad (\mbox{root})
\end{array}
\right.
$$
Similarly, in the case of $H=(e_1-e_2)\pm (e_2-e_i)\;(3\leq i)$,
we can put
$$
H = (e_1-e_2)\pm (e_2-e_3)
= \left\{\begin{array}{l}
e_1 - 2e_2 + e_3 \sim e_1 + e_2 - 2e_3 \\
e_1 - e_3 \quad (\mbox{root})
\end{array}
\right.
$$
Here, for $H_1, H_2 \in \mathfrak a$, we express $H_1 \sim H_2$
when $H_1$ can be translated to $H_2$ by some element of $K$.
In other words, $H_1$ is equivalent to $H_2$
under the action of the Weyl group.

In the case of $H=(e_1-e_2)\pm (e_i-e_j)\; (3\leq i<j)$,
there exists an element of the Weyl group 
which fixes $e_1, e_2$ and translates $e_i$ to $e_3$ and $e_j$ to $e_4$.
Therefore we can put
$$
H = (e_1-e_2) \pm (e_3-e_4)
= \left\{\begin{array}{l} e_1+e_3-e_2-e_4 \\
e_1+e_4-e_2-e_3 \end{array} \right.
$$
By the action of the Weyl group,
these vectors are equivalent to $e_1+e_2-e_3-e_4$.
Consequently, it suffices to consider orbits through
$$
H = 2e_1-e_2-e_3,\ e_1+e_2-2e_3,\ e_1+e_2-e_3-e_4,
$$
which have a possibility to be austere.

In the case of $H=2e_1-e_2-e_3$,
the only possibility to be the form
$H=\mbox{(a positive root)}\pm \mbox{(a positive root)}$ is
$H=(e_1-e_2)+(e_1-e_3)$.
Thus, from Lemma~\ref{lem:A}, the set $R_+-R_+^\Delta$ must be
$$
R_+-R_+^\Delta =\{e_1-e_2,\; e_1-e_3\}.
$$
When $l\geq 3$, since $\langle e_3-e_4, H\rangle\not= 0$,
we have $e_3-e_4\in R_+-R_+^\Delta$.
This is a contradiction.
Hence $l=2$ and then $\mathrm{Ad}(K)H$ is austere in $S$.
Similarly, the orbit through $H = e_1+e_2-2e_3$ is also austere.

In the case of $H=e_1+e_2-e_3-e_4$, possibilities to be the form
$H=\mbox{(a positive root)}\pm \mbox{(a positive root)}$ are
$$
H=(e_1-e_3)+(e_2-e_4)=(e_1-e_4)+(e_2-e_3).
$$
Thus $R_+-R_+^\Delta$ must satisfy
$$
R_+-R_+^\Delta\subset\{e_1-e_3,\; e_2-e_4,\; e_1-e_4,\; e_2-e_3\}
$$
When $l\geq 4$, since $\langle e_4-e_5,H\rangle\not= 0$,
we have $e_4-e_5\in R_+-R_+^\Delta$.
This is a contradiction.
Hence $l=3$, and then $\mathrm{Ad}(K)H$ is austere in $S$.
\end{proof}

\begin{pro} \label{pro:austere-D}
In the case where $R$ is of type $D_l$,
an austere orbit is one of the following
except orbits through a restricted root vector:
\begin{enumerate}
\item the orbit through $H=e_1$,
\item when $l=4$, the orbit through $H=e_1+e_2+e_3+e_4,\; e_1+e_2+e_3-e_4$.
\end{enumerate}
\end{pro}

\begin{proof}
In the case of $R=D_l$, 
$R_+$ is given by $R_+ = \{ e_i \pm e_j \mid i<j \}$.
Since all restricted roots have constant multiplicities,
the condition (\ref{eqn:multiplicity})
of Lemma~\ref{lem:A} is always satisfied.
It is easy to see that the orbit through $e_1$ (or its scalar multiple)
is austere.
Therefore we consider other orbits.
From Lemma~\ref{lem:A},
we can assume $H=\pm\mbox{(a positive root)}\pm\mbox{(a positive root)}$.
Since any root can be translated to $e_1+e_2$ by the action of the Weyl group,
we can assume $H=(e_1+e_2)\pm\mbox{(a positive root)}$.
Furthermore any root can be translated to one of
$$
e_1\pm e_2,\; e_1+e_3,\; e_2+e_4,\; e_3+e_4,\; e_3-e_4
$$
by the action of elements of the Weyl group which fix $e_1, e_2$.
Therefore $H$ is one of the following:
\begin{eqnarray*}
H &=& (e_1+e_2) \pm (e_1-e_2) = 2e_1,\ 2e_2 \sim 2e_1, \\
H &=& (e_1+e_2) \pm (e_1+e_3)
= \left\{ \begin{array}{l} 2e_1+e_2+e_3, \\
e_2-e_3 \ (\mbox{root}), \end{array} \right. \\
H &=& (e_1+e_2) \pm (e_2+e_4)
= \left\{ \begin{array}{l} 2e_2+e_1+e_4 \sim 2e_1+e_2+e_3, \\
e_1-e_4 \ (\mbox{root}), \end{array} \right. \\
H &=& (e_1+e_2) \pm (e_3+e_4) \sim e_1+e_2+e_3+e_4, \\
H &=& (e_1+e_2) \pm (e_3-e_4) \sim e_1+e_2+e_3-e_4.
\end{eqnarray*}
Consequently, it suffices to consider orbits through
$$
H=2e_1+e_2+e_3,\; e_1+e_2+e_3+e_4,\; e_1+e_2+e_3-e_4
$$
which have a possibility to be austere.

In the case of $H=2e_1+e_2+e_3$,
the only possibility to be the form
$H=\mbox{(a positive root)}\pm \mbox{(a positive root)}$ is
$$
H=(e_1+e_2)+(e_1+e_3).
$$
Thus $R_+-R_+^\Delta$ must be
$$
R_+-R_+^\Delta =\{e_1+e_2,\; e_1+e_3\}.
$$
Since $\langle e_1-e_2, H\rangle \neq 0$,
we have $e_1-e_2\in R_+-R_+^\Delta$.
This is a contradiction. Hence this orbit is not austere.

In the case of $H=e_1+e_2+e_3+e_4$,
possibilities of the form
$H=\mbox{(a positive root)}\pm \mbox{(a positive root)}$ are
$$
H = (e_1+e_2)+(e_3+e_4) = (e_1+e_3)+(e_2+e_4) = (e_1+e_4)+(e_2+e_3).
$$
Thus $R_+ - R_+^\Delta$ must satisfy
$$
R_+ - R_+^\Delta \subset
\{e_1+e_2,\; e_3+e_4,\; e_1+e_3,\; e_2+e_4,\; e_1+e_4,\; e_2+e_3\}.
$$
When $l\geq 5$, since $\langle e_4+e_5,H\rangle\not= 0$,
we have $e_4+e_5\in R_+-R_+^\Delta$.
This is a contradiction.
Hence $l=4$, and then the orbit $\mathrm{Ad}(K)H$ is austere in $S$.
In the case of $H=e_1+e_2+e_3-e_4$,
similarly we have $l=4$, and then $\mathrm{Ad}(K)H$ is austere.
\end{proof}

\begin{pro} \label{pro:austere-BC}
In the case where $R$ is of type $B_l, C_l$ or $BC_l$,
an austere orbit except orbits through a restricted root vector
is the following:

When $R=B_2$ where the multiplicities of the restricted roots
are constant, the orbit through
$$
H=e_1+\frac{e_1+e_2}{\sqrt{2}}
$$
is austere. This orbit is a principal orbit.
\end{pro}

\begin{rem} \label{rem:austere-BC} \rm
In the case of $R=B_2$,
there exist two singular orbits and these are not isometric.
Hence from Proposition
\ref{pro:principal orbit of cohomogeneity one actions},
a principal austere orbit in Proposition \ref{pro:austere-BC}
is not a weakly reflective submanifold.
\end{rem}

\begin{proof}
First we consider the case of $R=B_l$,
where $R_+=\{e_i,\; e_i\pm e_j \; | \; i<j \}$.
From Lemma~\ref{lem:A}, we can assume
$$
H = \frac{\alpha}{\| \alpha \|} \pm \frac{\beta}{\| \beta \|}\quad 
(\alpha ,\beta\in R_+).
$$

i)\;
When $\alpha$ and $\beta$ are both short roots,
we can put $H=\alpha\pm\beta$.
Furthermore, since any short root $\alpha$ can be translated to $e_1$
by the action of the Weyl group,
we can assume $H=e_1\pm\beta$.
If $\beta =e_1$, then $H=2e_1$ and this is equivalent to
the orbit through a root vector.
If $\beta =e_j\;(j\geq 2)$, then $H=e_1\pm e_j$ is a root vector.

ii)\:
When $\alpha$ and $\beta$ are both long roots,
we can put $H=\alpha\pm\beta$.
Since any long root $\alpha$ can be translated to $e_1+e_2$
by the action of the Weyl group,
we can assume $H=(e_1+e_2)\pm\beta$.
Furthermore $\beta$ can be translated to one of
$$
\beta = e_1 \pm e_2,\; e_1+e_3,\; e_2+e_3,\; e_3+e_4
$$
by the action of elements of the Weyl group
which fix $e_1$ and $e_2$.

In the case of $\beta = e_1 \pm e_2$, $H$ is equivalent to a root or zero vector.

In the case of $\beta = e_1+e_3$,
$$
H = (e_1+e_2) \pm (e_1+e_3) = \left\{ \begin{array}{l}
2e_1+e_2+e_3, \\ e_2-e_3 \quad (\mbox{root}).
\end{array} \right.
$$
When $H = 2e_1+e_2+e_3$, the only possibility to be the form
$H=\mbox{(a positive root)}\pm \mbox{(a positive root)}$ is 
$H = (e_1+e_2) + (e_1+e_3)$.
Thus $R_+-R_+^\Delta$ must be
$R_+-R_+^\Delta = \{ e_1+e_2,\; e_1+e_3 \}$.
On the other hand, since $\langle e_1, H \rangle \not= 0$,
we have $e_1 \in R_+-R_+^\Delta$.
This is a contradiction.
Hence this orbit is not austere.

In the case of $\beta =e_2+e_3$,
$$
H = (e_1+e_2) \pm (e_2+e_3)
= \left\{\begin{array}{l}
e_1 + 2e_2 + e_3\sim 2e_1 + e_2 + e_3, \\
e_1 - e_3 \quad (\mbox{root}).
\end{array}
\right.
$$

In the case of $\beta =e_3+e_4$,
$$
H = (e_1 + e_2) \pm (e_3 + e_4)
= \left\{\begin{array}{l}
e_1 + e_2 + e_3 + e_4, \\
e_1 + e_2 - e_3 - e_4 \sim e_1 + e_2 + e_3 + e_4.
\end{array}
\right.
$$
In this case,
$$
H=(e_1+e_2)+(e_3+e_4)=(e_1+e_3)+(e_2+e_4)=(e_1+e_4)+(e_2+e_3)
$$
are possibilities to be the form
$H = \alpha \pm \beta$.
Thus $R_+ - R_+^\Delta$ must satisfy
$$
R_+ - R_+^\Delta \subset
\{e_1+e_2,\; e_3+e_4,\; e_1+e_3,\; e_2+e_4,\; e_1+e_4,\; e_2+e_3\}.
$$
On the other hand,
since $\langle e_1,H\rangle\not=0$,
we have $e_1\in R_+ - R_+^\Delta$.
This is a contradiction.
Hence this orbit is not austere.

iii)\,
When $\alpha$ is a short root and $\beta$ is a long root,
we can assume $\alpha = e_1$ and
$$
H=e_1\pm\frac{\beta}{\sqrt{2}} \quad \mbox{where} \quad \beta =e_1+e_2,\; e_2+e_3.
$$

In the case of $H = e_1 + \frac{e_1+e_2}{\sqrt 2}$,
if $l \geq 3$, then $e_3 \in R_+ - R_+^\Delta$.
On the other hand, there is no $\mu \in R_+$ such that
$$
H = n \left(e_3 \pm \frac{\mu}{\| \mu \|}\right).
$$
Thus we have $l=2$.
In this case 
$R_+ - R_+^\Delta = \{e_1,\; e_2,\; e_1 + e_2,\; e_1 - e_2\}$.
If we define $f:R_+-R_+^\Delta\rightarrow R_+-R_+^\Delta$ by
$$
f(e_1) = e_1 + e_2, \qquad f(e_2) = e_1 - e_2,
$$
then $H$ satisfies the condition (\ref{eqn:1}) of Lemma~\ref{lem:A}.
Hence this orbit is austere if the multiplicities of the restricted roots
are constant.

In the case of $H=e_1-\frac{e_1+e_2}{\sqrt{2}}$, we can express $H$ as
$$
H=-\frac{1}{\sqrt{2}+1}\left\{
\left(1+\frac{1}{\sqrt{2}}\right)e_2-\frac{1}{\sqrt{2}}e_1
\right\}.
$$
Permuting $e_1$ and $e_2$ by the action of the Weyl group
and replacing $e_2\mapsto -e_2$, we have that this orbit is equivalent to
the orbit through
$$
H=e_1+\frac{e_1+e_2}{\sqrt{2}}.
$$

In the case of $\beta =e_2+e_3$,
$$
H=e_1\pm\frac{e_2+e_3}{\sqrt{2}} \sim e_1+\frac{e_2+e_3}{\sqrt{2}}.
$$
In this case $e_3 \in R_+ - R_+^\Delta$.
On the other hand, there is no $\mu \in R_+$ such that
$$
H = n \left(e_3 \pm \frac{\mu}{\| \mu \|}\right).
$$
Hence this orbit is not austere.

Second we consider the case of $R=C_l$,
where $R_+=\{2e_i,\; e_i\pm e_j \; | \; i<j \}$.
For this purpose we shall use the dual mapping 
and transfer the result of the case $R=B_l$ by the dual mapping.
A mapping
$$
\frak a-\{0\}\rightarrow\frak a-\{0\};
H\mapsto H^*=\frac{2H}{\langle H,H\rangle}
$$
is called a dual mapping.
This maps a root system to a root system,
more precisely, a long root is moved to a short root and a short root
is moved to a long root.
Root systems of type $B_l$ and $C_l$ are dual
by this mapping, and other irreducible root systems are self-dual.
If there exists $f$ which satisfies (\ref{eqn:1}) for $H \in \mathfrak a$,
then there exists $f^*$ which satisfies (\ref{eqn:1}) for $H^*$.

In the above discussion, in the case of $R=B_l\ (l \geq 3)$,
we showed that there are no austere orbits
except orbits through a restricted root vector.
Thus we also have that there are no austere orbits
except orbits through a restricted root vector in the case of $R=C_l\ (l \geq 3)$.
When $l=2$, $C_2=B_2$.

Finally we shall consider the case of $R=BC_l$,
where $R_+=\{e_i,\; 2e_i,\; e_i\pm e_j \; | \; i<j \}$.
From Lemma~\ref{lem:A}, we can put
$$
H = \frac{\alpha}{\| \alpha \|} \pm \frac{\beta}{\| \beta \|}\quad
(\alpha ,\beta\in R_+).
$$
If $\| \alpha \| = \| \beta \|$, then we can put $H = \alpha \pm \beta$.
When $\alpha$ and $\beta$ are both short roots or both long roots,
$H$ is a scalar multiple of a root vector.
When $\alpha$ and $\beta$ are both middle roots,
we can assume $\alpha = e_1+e_2$ and
$$
H = (e_1+e_2) \pm \beta \quad (\beta = e_1+e_3,\; e_2+e_3,\; e_3+e_4).
$$
By the action of the Weyl group, these are equivalent to
$$
H=2e_1+e_2+e_3, \quad H=e_1+e_2+e_3+e_4
$$
or a restricted root vector.
In the case of $H=2e_1+e_2+e_3$, we have $l=3$ and
$$
R_+ - R_+^\Delta = \{e_1,\; e_2,\; e_3,\; 2e_1,\; 2e_2,\; 2e_3,\;
e_1 \pm e_2,\; e_1 \pm e_3,\; e_2+e_3\}.
$$
Since there is no $f$ which satisfies (\ref{eqn:1}),
this orbit is not austere.
In the case of $H = e_1+e_2+e_3+e_4$, we have $l=4$ and
\begin{eqnarray*}
R_+ - R_+^\Delta &=& \{e_i,2e_i\;|\;1\leq i\leq 4\} \\
& & \cup \{e_1+e_2,\; e_1+e_3,\; e_1+e_4,\; e_2+e_3,\; e_2+e_4,\; e_3+e_4\}.
\end{eqnarray*}
Since there is no $f$ which satisfies (\ref{eqn:1}),
this orbit is not austere.

It remains the case where $\| \alpha \| < \| \beta \|$.
When $\alpha$ is a short root and $\beta$ is a long root,
$H$ is a scalar multiple of a root vector.
By the dual mapping, we can identify two cases,
where $\alpha$ is a short root and $\beta$ is a middle root,
and where $\alpha$ is a middle root and $\beta$ is a long root.
Therefore we shall discuss the former.
In this case, we can assume $\alpha =e_1$ and
$$
H=e_1\pm\frac{\beta}{\sqrt{2}}\quad \mbox{where} \quad \beta =e_1+e_2,\; e_2+e_3.
$$
Similarly with the case of the restricted root system of type $B$,
$$
H=e_1+\frac{e_1+e_2}{\sqrt{2}}, \qquad
H=e_1+\frac{e_2+e_3}{\sqrt{2}}
$$
have a possibility to be austere.
When $H=e_1+\frac{e_2+e_3}{\sqrt{2}}$, we have $l=3$ and 
$$
R_+ - R_+^\Delta = \{e_1,\; e_2,\; e_3,\; 2e_1,\; 2e_2,\; 2e_3,\;
e_1 \pm e_2,\; e_1 \pm e_3,\; e_2+e_3\}.
$$
Since there is no $f$ which satisfies (\ref{eqn:1}),
this orbit is not austere.
When $H=e_1+\frac{e_1+e_2}{\sqrt{2}}$,
the orbit has a possibility to be austere if $l=2$.
In this case, the orbit is a principal orbit.
This orbit is austere if
the sum of the multiplicities of long roots and short roots
coincides with the multiplicity of middle roots.
From the classification of symmetric pairs,
there does not exist such a symmetric pair.
\end{proof}

\begin{pro} \label{pro:austere-G}
In the case where $R$ is of type $G_2$,
the orbit through
$$
H=\alpha_1+\frac{\alpha_2}{\sqrt{3}}
$$
is the only austere orbit except orbits through a restricted root vector.
This orbit is a principal orbit.
\end{pro}

\begin{rem} \rm
By the same discussion in Remark \ref{rem:austere-BC},
this principal austere orbit is not a weakly reflective submanifold
from Proposition~\ref{pro:principal orbit of cohomogeneity one actions}.
\end{rem}

\begin{proof}
The fundamental system $F$ of the restricted root system of type $G_2$
is given by
$F = \{\alpha_1=e_1-e_2,\; \alpha_2=-2e_1+e_2+e_3\}$
and the set $R_+$ of positive roots is
$$
R_+ = F \cup \left\{
\begin{array}{l} \alpha_1+\alpha_2=-e_1+e_3,\; 2\alpha_1+\alpha_2=-e_2+e_3, \\
3\alpha_1+\alpha_2=e_1-2e_2+e_3,\; 3\alpha_1+2\alpha_2=-e_1-e_2+2e_3 \end{array}
\right\}.
$$
In the case of $G_2$, since all restricted roots have constant multiplicities,
the condition (\ref{eqn:multiplicity})
of Lemma~\ref{lem:A} is always satisfied.
From Lemma~\ref{lem:A}, we can put
$$
H = \frac{\alpha}{\| \alpha \|} \pm \frac{\beta}{\| \beta \|} \quad
(\alpha ,\beta \in R_+).
$$

When $\alpha$ and $\beta$ are both short roots,
we can put $H=\alpha\pm\beta\; (\alpha\not=\beta )$.
Furthermore, since any short root $\alpha$ can be translated to $\alpha_1$
by the action of the Weyl group, we can assume
$$
H=\alpha_1\pm\beta \quad \mbox{where} \quad
\beta=\alpha_1+\alpha_2,\; 2\alpha_1+\alpha_2.
$$
In the case of $\beta =\alpha_1+\alpha_2$,
$$
H=\alpha_1\pm (\alpha_1+\alpha_2)=2\alpha_1+\alpha_2,\; -\alpha_2.
$$
Then $H$ is a root vector.
In the case of $\beta=2\alpha_1+\alpha_2$,
$$
H=\alpha_1\pm (2\alpha_1+\alpha_2)=3\alpha_1+\alpha_2,\; -\alpha_1-\alpha_2.
$$
Then $H$ is a root vector.

When $\alpha$ and $\beta$ are both long roots,
we can put
$H=\alpha\pm\beta\; (\alpha\not=\beta )$.
Since any long root $\alpha$ can be translated to $\alpha_2$
by the action of the Weyl group, we can assume
$$
H=\alpha_2\pm\beta \quad \mbox{where} \quad
\beta=3\alpha_1+\alpha_2,\; 3\alpha_1+2\alpha_2.
$$
In the case of $\beta=3\alpha_1+\alpha_2$,
$$
H=\alpha_2\pm (3\alpha_1+\alpha_2)=-3\alpha_1,\; 3\alpha_1+2\alpha_2.
$$
Then $H$ is a scalar multiple of a root vector.
In the case of $\beta=3\alpha_1+2\alpha_2$,
$$
H=\alpha_2\pm (3\alpha_1+2\alpha_2)=-3\alpha_1-\alpha_2,\; 3(\alpha_1+\alpha_2).
$$
Then $H$ is a scalar multiple of a root vector.

When $\alpha$ is a short root and $\beta$ is a long root,
we can assume $\alpha =\alpha_1$ and
$$
H=\alpha_1\pm \frac{\beta}{\sqrt{3}}\quad \mbox{where} \quad
\beta =\alpha_2,\; 3\alpha_1+\alpha_2,\; 3\alpha_1+2\alpha_2.
$$
We note that the orbit though $H$ is a principal orbit.
In these cases, $H$ is equivalent to a scalar multiple of 
$\alpha_1+\frac{\alpha_2}{\sqrt{3}}$ by the action of the Weyl group.
In the case of $H=\alpha_1+\frac{\alpha_2}{\sqrt{3}}$,
$$
H=\frac{1}{\sqrt{3}+1}\left(
\alpha_1+\alpha_2+\frac{3\alpha_1+\alpha_2}{\sqrt{3}}
\right)
=\frac{1}{\sqrt{3}+2}\left(
2\alpha_1+\alpha_2+\frac{3\alpha_1+2\alpha_2}{\sqrt{3}}
\right).
$$
Thus from Lemma~\ref{lem:A} this orbit is austere.
This completes the proof.
\end{proof}

\begin{pro}
In the case of $R=F_4$, there are no austere orbits
except orbits through a restricted root vector.
\end{pro}

\begin{proof}
In this case $R_+$ is given by
\begin{eqnarray*}
R_+ = \{ e_i \}_{1 \leq i \leq 4} \cup \{ e_i \pm e_j \}_{1 \leq i < j \leq 4}
\cup \left\{ \frac{1}{2}(e_1 \pm e_2 \pm e_3 \pm e_4) \right\}
\end{eqnarray*}
From Lemma~\ref{lem:A}, we can assume
$$
H = \frac{\alpha}{\| \alpha \|} \pm \frac{\beta}{\| \beta \|}\quad 
(\alpha ,\beta \in R_+).
$$
When $\alpha$ and $\beta$ are both short roots,
we can put $H = \alpha \pm \beta$.
In this case, since any short root can be translated to $e_1$
by the action of the Weyl group,
we can put
$$
H = e_1 \pm \beta \quad \mbox{where} \quad
\beta =e_i\ (i\geq 2),\; \frac{1}{2}(e_1\pm e_2\pm e_3 \pm e_4).
$$
In the case of $\beta=e_i$, $H$ is a root vector.
In the case of
$\beta = \frac{1}{2}(e_1\pm e_2\pm e_3\pm e_4)$.
$$
H = e_1 \pm \frac{1}{2}(e_1 \pm e_2 \pm e_3 \pm e_4)
= \left\{ \begin{array}{l}
e_1 + \frac{1}{2}(e_1 \pm e_2 \pm e_3 \pm e_4), \\
\frac{1}{2}(e_1 \pm e_2 \pm e_3 \pm e_4) \quad \mbox{(a root)}.
\end{array} \right.
$$
Thus we consider the case of 
$H = e_1 + \frac{1}{2}(e_1 \pm e_2 \pm e_3 \pm e_4)$.
In this case, $\langle e_4, H \rangle \not= 0$,
however, there does not exist $\beta \in R_+$ such that
$$
H = n \left( e_4 \pm \frac{\beta}{\| \beta \|} \right).
$$
Hence from Lemma~\ref{lem:A} this orbit is not austere.

When $\alpha$ and $\beta$ are both long root,
we can assume $H=e_1+e_2\pm e_i\pm e_j\;(i<j)$.
Moreover, we exclude $H$ which is a scalar multiple of a root vector.
Then 
$$
H = 2e_1 + e_2 \pm e_i,\ e_1 + 2e_2 \pm e_i\ (i=3,4).
$$
The reflection $s_{e_3-e_4}$ permutes $e_3$ and $e_4$,
and fixes $e_1,e_2$.
Therefore we can put
$$
H = 2e_1 + e_2 \pm e_3.
$$
In this case, $\langle e_3, H \rangle \not= 0$,
however, there does not exist $\beta \in R_+$ such that
$$
H = n \left( e_3 \pm \frac{\beta}{\| \beta \|} \right).
$$
Hence this orbit is not austere.

When $\alpha$ is a short root and $\beta$ is a long root,
we can put
$$
H=e_1+\frac{\pm e_i\pm e_j}{\sqrt{2}}\quad (i<j).
$$
Moreover, by the action of the Weyl group we can assume
$$
H=e_1+\frac{\pm e_1+e_2}{\sqrt{2}}\quad\mbox{or}\quad 
H=e_1+\frac{e_2+e_3}{\sqrt{2}}.
$$
In the case of $H=e_1+\frac{e_2+e_3}{\sqrt{2}}$,
we have
$\langle e_3, H \rangle \not = 0$.
However, there does not exist $\beta \in R_+$ such that
$$
H = n \left( e_3 \pm \frac{\beta}{\| \beta \|} \right).
$$
Hence this orbit is not austere.
In the case of $H = e_1 + \frac{\pm e_1+e_2}{\sqrt{2}}$,
we have
$\langle e_2+e_3,H\rangle\not=0$.
However, there does not exist $\beta \in R_+$ such that
$$
H=n\left(\frac{e_2+e_3}{\sqrt{2}}+\frac{\beta}{\| \beta \|}\right).
$$
Hence this orbit is not austere.
\end{proof}

\begin{pro}
In the case of $R=E_8$, there are no austere orbits
except the orbits through a restricted root vector.
\end{pro}

\begin{proof}
In the case of $R=E_8$, $R_+$ is given by
$$
R_+=\{\pm e_i+e_j\;|\;1\leq i<j\leq 8\}\cup
\left\{\left.\frac{1}{2}(e_8+\sum_{i=1}^7(-1)^{\nu (i)}e_i)\;\right|\;
\sum_{i=1}^7\nu (i)\in 2\mbox{\boldmath$Z$}\right\}
$$
Since all restricted roots have constant multiplicities,
the condition (\ref{eqn:multiplicity})
of Lemma~\ref{lem:A} is always satisfied.
From Lemma~\ref{lem:A},
we can put $H=e_1+e_2+\beta$ where
$$
\beta =\left\{
\begin{array}{ll}
\pm e_1\pm e_2, & \\
\pm e_1\pm e_i & (3\leq i\leq 8),\\
\pm e_2\pm e_i & (3\leq i\leq 8),\\
\pm e_i\pm e_j & (3\leq i<j\leq 8),\\
\frac{1}{2}\sum_{i=1}^{\nu (i)}(-1)^{\nu (i)}e_i & 
(\sum_{i=1}^8 \nu(i) \in 2 \mbox{\boldmath$Z$}).
\end{array}
\right.
$$

i) In the case of $\beta =\pm e_1\pm e_2$,
$$
H=\left\{
\begin{array}{l}
e_1+e_2+e_1+e_2=2(e_1+e_2)\quad\mbox{(twice of a root)},\\
e_1+e_2+e_1-e_2=2e_1,\\
e_1+e_2-e_1+e_2=2e_2,\\
e_1+e_2-e_1-e_2=0.
\end{array}
\right.
$$
When $H=2e_1$,
$$
R-R^\Delta 
= \{\pm e_1\pm e_j\;|\;2\leq j\leq 8\}\cup 
\left\{\left.\pm\frac{1}{2}(e_1+\sum_{i=2}^8(-1)^{\nu (i)}e_i)\;\right|\;
\sum_{i=2}^8 \nu (i)\in 2\mbox{\boldmath$Z$}\right\}
$$
Then there does not exist $\beta \in R_+$ such that
$H = n(\frac{1}{2}\sum_{i=1}^8 e_i\pm\beta)$.
Hence this orbit is not austere.
When $H=2e_2$,
since the reflection $s_{e_1-e_2}$ permutes $e_1$ and $e_2$,
this orbit is equivalent to the orbit through $H=2e_1$.
Hence this orbit is not austere.

ii) In the case of $\beta=\pm e_1\pm e_i\;(3\leq i\leq 8)$,
$$
H = \left\{ \begin{array}{l}
2e_1+e_2\pm e_i, \\
e_2\pm e_i\quad\mbox{(a root)}.
\end{array}
\right.
$$
The reflection $s_{e_3-e_i}\;(i\geq 4)$ fixes $e_1,e_2$ and permutes
$e_3$ and $e_i$.
Thus the orbit through $H=2e_1+e_2\pm e_i$ is equivalent to
the orbit through $H=2e_1+e_2\pm e_3$.
Then $\langle e_1+e_4,H\rangle\not= 0$, however,
there does not exist $\beta\in R_+$ such that $H=n(e_1+e_4\pm\beta )$.
Hence this orbit is not austere.

iii) In the case of $\beta=\pm e_2\pm e_i\; (3 \leq i \leq 8)$,
$$
H = \left\{ \begin{array}{l}
e_1+2e_2\pm e_i, \\
e_1\pm e_i\quad\mbox{(a root)}.
\end{array}
\right.
$$
By the action of the Weyl group,
the orbit through $H=e_1+2e_2\pm e_i$ is equivalent to the orbit
through $H=2e_1+e_2\pm e_i$.
Hence this orbit is not austere.

iv) In the case of $\beta =\pm e_i\pm e_j\; (3 \leq i < j \leq 8)$,
we can assume
$$
H=e_1+e_2\pm e_3\pm e_4.
$$
Then $\langle e_1+e_5,H\rangle\not= 0$, however,
there does not exist $\beta \in R_+$ such that $H=n(e_1+e_5\pm\beta )$.
Thus this orbit is not austere.

v) v-1) In the case of
$$
\beta =\frac{1}{2}(-e_1-e_2+\sum_{i=3}^8 (-1)^{\nu (i)}e_i),\quad 
\sum_{i=3}^8 \nu (i)\in 2\mbox{\boldmath$Z$},
$$
then
$$
H=\frac{1}{2}(e_1+e_2+\sum_{i=3}^8 (-1)^{\nu (i)}e_i)\quad\mbox{(a root)}.
$$
v-2) In the case of
$$
\beta =\frac{1}{2}(-e_1+e_2+\sum_{i=3}^8 (-1)^{\nu (i)}e_i),\quad 
\sum_{i=3}^8 \nu (i)\in 2\mbox{\boldmath$Z$}+1,
$$
then
$$
H=\frac{1}{2}(e_1+3e_2+\sum_{i=3}^8 (-1)^{\nu (i)}e_i)\sim 
\frac{1}{2}(3e_1+e_2+\sum_{i=3}^8 (-1)^{\nu (i)}e_i).
$$
v-3) In the case of
$$
\beta =\frac{1}{2}(e_1-e_2+\sum_{i=3}^8 (-1)^{\nu (i)}e_i),\quad 
\sum_{i=3}^8 \nu (i)\in 2\mbox{\boldmath$Z$}+1,
$$
then
$$
H=\frac{1}{2}(3e_1+e_2+\sum_{i=3}^8 (-1)^{\nu (i)}e_i).
$$
In this case, there exists an $i$ such that $\nu(i) = 1$.
Permuting $e_i$ and $e_3$ by the action of the Weyl group,
we have
$$
H=\frac{1}{2}(3e_1+e_2-e_3+\sum_{i=4}^8 (-1)^{\nu (i)}e_i),\quad 
\sum_{i=4}^8 \nu (i)\in 2\mbox{\boldmath$Z$}.
$$
Then $\langle e_2-e_3,H\rangle\not= 0$, however,
there does not exist $\beta\in R_+$ such that $H=n(e_2-e_3\pm\beta )$.
Thus this orbit is not austere.

vi) In the case of
$$
\beta =\frac{1}{2}(e_1+e_2+\sum_{i=3}^8 (-1)^{\nu (i)}e_i),\quad 
\sum_{i=3}^8 \nu (i)\in 2\mbox{\boldmath$Z$},
$$
then
$$
H=\frac{1}{2}(3e_1+3e_2+\sum_{i=3}^8 (-1)^{\nu (i)}e_i).
$$
In this case $\langle e_1+e_3,H\rangle\not=0$, however,
there does not exist $\beta\in R_+$ such that $H = n(e_1+e_3\pm\beta)$.
Thus this orbit is not austere.
\end{proof}

\begin{pro}
In the case of $R=E_7$, there are no austere orbits
except the orbits through a restricted root vector.
\end{pro}

\begin{proof}
In the case of $R=E_7$, all restricted roots have constant multiplicities.
Thus the condition (\ref{eqn:multiplicity})
of Lemma~\ref{lem:A} is always satisfied.
$\frak a=\{\sum_{i=1}^8\xi_ie_i\;|\;\xi_8=-\xi_7\}$
and
\begin{eqnarray*}
R_+&=&\{\pm e_i+e_j\;|\;1\leq i<j\leq 6\}\cup\{e_8-e_7\}\\
& &\cup\left\{\left.\frac{1}{2}(e_8-e_7+\sum_{i=1}^6 (-1)^{\nu (i)}e_i)
\;\right|\;\sum_{i=1}^6 \nu (i)\in 2\mbox{\boldmath$Z$}+1\right\}
\end{eqnarray*}
From Lemma~\ref{lem:A},
we can assume $H=e_7-e_8+\beta$ where
$$
\beta =\left\{
\begin{array}{l}
\pm e_i\pm e_j\quad (1\leq i<j\leq 6),\\
\pm\frac{1}{2}(e_7-e_8+\sum_{i=1}^6(-1)^{\nu (i)}e_i).
\end{array}
\right.
$$

In the case of $\beta = \pm e_i\pm e_j\; (1\leq i<j\leq 6)$,
we take $k$ such that $1\leq k\leq 6$ and $k\not= i,j$.
Then $e_j+e_k \in R_+$ and $\langle H, e_j+e_k \rangle \not= 0$,
however, there does not exist $\alpha \in R$ such that $H=n(e_j+e_k+\alpha)$.
Thus this orbit is not austere.

In the case of
$\beta = \pm\frac{1}{2}(e_7-e_8+\sum_{i=1}^6(-1)^{\nu (i)}e_i)$
then
$$
H = \left\{ \begin{array}{l}
\frac{1}{2}(3e_7-3e_8+\sum_{i=1}^6(-1)^{\nu (i)}e_i),\\
\frac{1}{2}(e_7-e_8-\sum_{i=1}^6(-1)^{\nu (i)}e_i)\quad\mbox{(a root)}.
\end{array}
\right. 
$$
Therefore it suffices to consider the case of
$H=\frac{1}{2}(3e_7-3e_8+\sum_{i=1}^6(-1)^{\nu (i)}e_i)$.
In this case, either $e_1+e_2$ or $e_1-e_2$ is an element of $R_+-R_+^\Delta$.
Denote this element by $\alpha$.
Then there does not exist $\beta \in R_+$ such that $H=n(\alpha\pm\beta)$.
Hence this orbit is not austere.
This completes the proof.
\end{proof}

\begin{pro}
In the case of $R=E_6$, there are no austere orbits
except orbits through a restricted root vector.
\end{pro}

\begin{proof}
In the case of $R=E_6$, all restricted roots have constant multiplicities.
Thus the condition (\ref{eqn:multiplicity})
of Lemma~\ref{lem:A} is always satisfied.
$\frak a=\{\sum_{i=1}^8\xi_ie_i\;|\;\xi_6=\xi_7=-\xi_8\}$
and
\begin{eqnarray*}
R_+ &=& \{\pm e_i + e_j\;|\;1\leq i<j\leq 5\}\\
& &\cup
\left\{\left.\frac{1}{2}(e_8-e_7-e_6+\sum_{i=1}^5 (-1)^{\nu (i)}e_i)
\;\right|\;\sum_{i=1}^5 \nu (i)\in 2\mbox{\boldmath$Z$}\right\}
\end{eqnarray*}
From Lemma~\ref{lem:A}, we can assume $H=e_1+e_2+\beta$ where
$$
\beta =\left\{
\begin{array}{l}
\pm (e_1-e_2),\\
\pm e_2\pm e_i\quad (3\leq i\leq 5),\\
\pm e_i\pm e_j\quad (3\leq i<j\leq 5),\\
\pm\frac{1}{2}(e_8-e_7-e_6+\sum_{i=1}^5 (-1)^{\nu (i)}e_i).
\end{array}
\right.
$$

i) In the case of $\beta =\pm (e_1-e_2)$, then $H=2e_1,\; 2e_2$.
For
$$
\alpha =\frac{1}{2}(e_8-e_7-e_6+e_1+e_2+e_3-e_4-e_5)\in R_+-R_+^\Delta
$$
there does not exist $\beta\in R_+$ such that $H=n(\alpha \pm \beta)$.
Thus this orbit is not austere.

ii) In the case of $\beta =\pm e_2\pm e_i\; (3\leq i\leq 5)$, then
$$
H = e_1 + e_2 \pm e_2 \pm e_i
= \left\{ \begin{array}{l} e_1+2e_2\pm e_i, \\ e_1 \pm e_i \quad (\mbox{root}).
\end{array} \right.
$$
Therefore it suffices to consider the case of $H = e_1+2e_2\pm e_i$. 
In this case, for $j$ with $3 \leq j \leq 5$ and $j \neq i$,
we have $e_1 + e_j \in R_+ - R_+^\Delta$.
However, there does not exist $\beta\in R_+$ such that $H=n(e_1+e_j \pm \beta)$.
Thus this orbit is not austere.

iii) In the case of $\beta =\pm e_i\pm e_j\; (3 \leq i < j \leq 5)$,
then $H=e_1+e_2\pm e_i\pm e_j$.
For $k$ with $3\leq k\leq 5$ and $k\not= i,j$,
we have $e_1+e_k \in R_+ - R_+^\Delta$.
However, there does not exist $\beta \in R_+$ such that $H=n(e_1+e_k \pm \beta)$.
Thus this orbit is not austere.

iv) In the case of
$$
\beta =\pm\frac{1}{2}(e_8-e_7-e_6+\sum_{i=1}^5 (-1)^{\nu (i)}e_i) \quad
\mbox{where} \quad \sum_{i=1}^5 \nu (i)\in 2\mbox{\boldmath$Z$},
$$
then
$$
H=e_1+e_2\pm\frac{1}{2}(e_8-e_7-e_6+\sum_{i=1}^5 (-1)^{\nu (i)}e_i).
$$
Here
$$
H=\left\{
\begin{array}{l}
\displaystyle{e_1+e_2+\frac{1}{2}
(e_8-e_7-e_6-e_1-e_2+\sum_{i=3}^5(-1)^{\nu (i)}e_i)},\\
\displaystyle{e_1+e_2-\frac{1}{2}
(e_8-e_7-e_6+e_1+e_2+\sum_{i=3}^5(-1)^{\nu (i)}e_i)}
\end{array}
\right.
$$
are root vectors.

In the case of
$$
H=\left\{
\begin{array}{l}
\displaystyle{e_1+e_2+\frac{1}{2}
(e_8-e_7-e_6+e_1-e_2+\sum_{i=3}^5(-1)^{\nu (i)}e_i)},\\
\displaystyle{e_1+e_2+\frac{1}{2}
(e_8-e_7-e_6-e_1+e_2+\sum_{i=3}^5(-1)^{\nu (i)}e_i)},\\

\displaystyle{e_1+e_2-\frac{1}{2}
(e_8-e_7-e_6+e_1-e_2+\sum_{i=3}^5(-1)^{\nu (i)}e_i)},\\
\displaystyle{e_1+e_2-\frac{1}{2}
(e_8-e_7-e_6-e_1+e_2+\sum_{i=3}^5(-1)^{\nu (i)}e_i)},
\end{array}
\right.
$$
we have $\langle e_1-e_2,H\rangle\not=0$.
However there does not exist $\beta\in R_+$ such that $H=n(e_1-e_2 \pm \beta)$.
Thus this orbit is not austere.

Finally, in the case of
$$
H=\left\{
\begin{array}{l}
\displaystyle{e_1+e_2+\frac{1}{2}
(e_8-e_7-e_6+e_1+e_2+\sum_{i=3}^5(-1)^{\nu (i)}e_i)}, \\
\displaystyle{e_1+e_2-\frac{1}{2}
(e_8-e_7-e_6-e_1-e_2+\sum_{i=3}^5(-1)^{\nu (i)}e_i)}
\end{array}
\right.
$$
we have $\langle e_1+e_3,H\rangle\not= 0$.
However there does not exist $\beta\in R_+$ such that $H=n(e_1+e_3 \pm \beta)$.
Thus this orbit is not austere.
This completes the proof.
\end{proof}

By discussions above, we completed the proof of
Theorem~\ref{thm:austere} and Theorem~\ref{thm:weakly reflective orbit}.

\section{Miscellaneous results}

In this section,
we shall concern with some results on weakly reflective submanifolds
besides orbits of $s$-representations.

\begin{pro}
Let $M_1$ and $M_2$ be weakly reflective submanifolds
in Riemannian manifolds $\tilde M_1$ and $\tilde M_2$, respectively.
Then $M_1 \times M_2$ is a weakly reflective submanifold
in $\tilde M_1 \times \tilde M_2$.
\end{pro}

\begin{proof}
We take $(x_1, x_2) \in M_1 \times M_2$ and a normal vector
$(\xi_1, \xi_2) \in T_{(x_1, x_2)}^\perp(M_1 \times M_2)$.
Let $\sigma_{\xi_1}$ and $\sigma_{\xi_2}$ be reflections of $M_1$ and $M_2$
in $\tilde M_1$ and $\tilde M_2$ with respect to $\xi_1$ and $\xi_2$, respectively.
Then $\sigma_{\xi_1}\times\sigma_{\xi_2}$ is an isometry of
$\tilde M_1 \times \tilde M_2$ and satisfies
\begin{eqnarray*}
&&
(\sigma_{\xi_1}\times\sigma_{\xi_2})(x_1, x_2) = (x_1, x_2), \\
&&
d(\sigma_{\xi_1}\times\sigma_{\xi_2})_{(x_1, x_2)}(\xi_1, \xi_2)
= -(\xi_1, \xi_2), \\
&&
(\sigma_{\xi_1}\times\sigma_{\xi_2})(M_1 \times M_2)
= M_1 \times M_2.
\end{eqnarray*}
Thus $M_1 \times M_2$ is a weakly reflective submanifold
in $\tilde M_1 \times \tilde M_2$.
\end{proof}

The following proposition states that the cone over a weakly reflective
submanifold in a sphere is also a weakly reflective submanifold in
a Euclidean space.

\begin{pro}
Let $M$ be a weakly reflective submanifold in a unit sphere $S^{n-1}(1)$.
Then the cone 
$$
C(M) = \{tx \mid t \in \mathbf R,\; t > 0,\; x \in M\}
$$
over $M$ is a weakly reflective submanifold
in a Euclidean space $\mathbf R^n$.
\end{pro}

\begin{proof}
Fix $x \in M$.
We note that for arbitrary $t \in \mathbf R\; (t > 0)$, we have
$$
T_{tx}^\perp(C(M)) = T_x^\perp M \subset T_xS^{n-1}(1).
$$
For $\xi \in T_x^\perp M$,
a reflection $\sigma_\xi$ of $M$ with respect to $\xi$
satisfies
$$
\sigma_\xi(x) = x, \qquad
(d\sigma_\xi)_x \xi = - \xi, \qquad
\sigma_\xi(M) = M.
$$
Since $\sigma_\xi$ is an isometry of $S^{n-1}(1)$,
it can be expressed as an orthogonal matrix.
Thus $\sigma_\xi$ acts on $\mathbf R^n$ and satisfies
$$
\sigma_\xi(tx) = t\sigma_\xi(x) = tx, \qquad
(d\sigma_\xi)_x \xi = - \xi.
$$
In addition,
for $x' \in M,\; t' \in \mathbf R,\; t' > 0$, we have
$$
\sigma_\xi(t'x') = t'\sigma_\xi(x') \in C(M).
$$
Therefore $\sigma_\xi(C(M)) = C(M)$.
Hence $C(M)$ is a weakly reflective submanifold in $\mathbf R^n$.
\end{proof}

Next we shall describe the relationship between 
weakly reflective submanifolds in an odd dimensional sphere and
in a complex projective space.

\begin{pro} \label{pro:wr-submfd in CP^n}
Let $M$ be a weakly reflective submanifold
in the hypersphere $S \subset \mathbf C^{n+1}$.
If $M$ is invariant under $U(1)$-actions and if a reflection of $M$
with respect to each normal vector is a unitary action,
then the image $P(M)$ of $M$
is a weakly reflective submanifold in $\mathbf CP^n$,
where $P$ is the natural projection $P : S \to \mathbf CP^n$.
\end{pro}

\begin{proof}
By the definition of $P$,
$T_x(U(1)x) = \ker dP_x$ at each point $x \in S$.
Thus $dP_x : T_x^\perp(U(1)x) \to T_x(\mathbf CP^n)$
is an isometric linear isomorphism.
Since $M$ is invariant under the $U(1)$-actions,
$P(M)$ is a submanifold in $\mathbf CP^n$.
Moreover, since $T_x(U(1)x) \subset T_xM$ for $x \in M$,
we have $T_x^\perp(U(1)x) \supset T_x^\perp M$.
Thus $dP_x : T_x^\perp M \to T_{p(x)}^\perp(P(M))$
also gives an isometric linear isomorphism.
Let $\sigma_\xi$ denote a reflection of $M$
with respect to $\xi \in T_x^\perp M$.
From the assumption,
$\sigma_\xi$ is a unitary transformation of $\mathbf C^{n+1}$.
Hence $\sigma_\xi$ induces an isometry of $\mathbf CP^n$.
Since $\sigma_\xi(x) = x$ and $\sigma_\xi(M) = M$,
we have $\sigma_\xi(P(x)) = P(x)$ and $\sigma_\xi(P(M)) = P(M)$.
In addition we have
\begin{eqnarray*}
d\sigma_\xi(dP_x(\xi))
&=& d(\sigma_\xi\circ P)_x(\xi)
= d(P\circ\sigma_\xi)_x(\xi) \\
&=& dP_x\sigma_\xi(\xi)
= dP_x(-\xi)
= - dP_x(\xi).
\end{eqnarray*}
Thus $\sigma_\xi$ is a reflection of $P(M)$ with respect to 
a normal vector $dP_x(\xi)$ at $P(x)$.
Hence $P(M)$ is a weakly reflective submanifold in $\mathbf CP^n$.
\end{proof}

\begin{cor}
An orbit of the $s$-representation
of an irreducible compact Hermitian symmetric pair
through a restricted root vector induces a weakly reflective submanifold in
a complex projective space.
\end{cor}

\begin{proof}
The center of the linear isotropy subgroup
of an irreducible compact Hermitian symmetric pair
is $U(1)$ (\cite{Helgason}).
Thus all orbits are invariant under $U(1)$.
Furthermore an orbit in the hypersphere $S$ through a restricted root vector
is a weakly reflective submanifold.
From the proofs of Lemma \ref{lem:orbits of s-representation}
and Proposition \ref{pro:orbit of root},
a reflection of an orbit through a restricted root vector
with respect to each normal vector is a unitary transformation.
Hence, from Proposition \ref{pro:wr-submfd in CP^n},
we have the conclusion.

\end{proof}


\end{document}